\documentclass[preprint,a4paper,12pt]{elsarticle}

\usepackage{amsmath}
\usepackage{amsfonts}
\usepackage{amssymb}
\usepackage{amsthm}

\usepackage{graphicx}
\graphicspath{{Figures/}}
\usepackage{subcaption}
\usepackage{multirow} % Table
\usepackage{float} % Position

\theoremstyle{definition}

\usepackage{xcolor}
\usepackage[colorlinks=true, allcolors=blue]{hyperref}

\newcommand{\bracketround}[1]{\left(#1\right)}
\newcommand{\bracketcurly}[1]{\left\{#1\right\}}
\newcommand{\bracketsquare}[1]{\left[#1\right]}

\newcommand{\norm}[1]{\left\lVert#1\right\rVert}

\journal{Applied Mathematics and Computation}

\begin{document}
	
	\begin{frontmatter}
		
		\title{A finite difference method with symmetry properties \\ for the high-dimensional Bratu equation}
		
		\author{Muhammad Luthfi Shahab}
		
		\author{Hadi Susanto\corref{cor1}}
		\cortext[cor1]{Corresponding author.}
		\ead{hadi.susanto@yandex.com}
		
		\author{Haralampos Hatzikirou}
		
		\address{
			%organization=
			{Department of Mathematics, Khalifa University of Science \& Technology},\\
			% addressline={Address One}, 
			%city=
			{Abu Dhabi},
			%postcode=
			{PO Box 127788}, 
			% state={State One},
			%country=
			{United Arab Emirates}
		}
		
		\begin{abstract}
			Solving the three-dimensional (3D) Bratu equation is highly challenging due to the presence of multiple and sharp solutions. Research on this equation began in the late 1990s, but there are no satisfactory results to date. 
			To address this issue, we introduce a symmetric finite difference method (SFDM) which embeds the symmetry properties of the solutions into a finite difference method (FDM). This SFDM is primarily used to obtain more accurate solutions and bifurcation diagrams for the 3D Bratu equation. Additionally, we propose modifying the Bratu equation by incorporating a new constraint that facilitates the construction of bifurcation diagrams and simplifies handling the turning points.
			The proposed method, combined with the use of sparse matrix representation, successfully solves the 3D Bratu equation on grids of up to $301^3$ points. The results demonstrate that SFDM outperforms all previously employed methods for the 3D Bratu equation. Furthermore, we provide bifurcation diagrams for the 1D, 2D, 4D, and 5D cases, and accurately identify the first turning points in all dimensions.
			All simulations indicate that the bifurcation diagrams of the Bratu equation on the cube domains closely resemble the well-established behavior on the ball domains described by Joseph and Lundgren \cite{joseph1973quasilinear}. Furthermore, when SFDM is applied to linear stability analysis, it yields the same largest real eigenvalue as the standard FDM despite having fewer equations and variables in the nonlinear system.
		\end{abstract}
		
		
	\end{frontmatter}

	\section{Introduction}
	
	%\subsection{Pengantar Singkat Bratu, Nonlinear PDEs, Numerical Solutions}
	The Bratu equation \cite{bratu1914equations} is a nonlinear elliptic partial differential equation (PDE) that arises in the study of steady-state solutions of the solid fuel ignition model \cite{jacobsen2002liouville}.  This equation is particularly intriguing due to the appearance of multiple solutions, some of which have nontrivial structures often referred to as sharp solutions. Having many applications, such as in combustion theory \cite{bebernes2013mathematical}, electrospinning process \cite{wan2004thermo}, chemical reactor theory \cite{boyd1986analytical}, electrostatics and plasma physics \cite{hichar2015application}, elastic nonlinear stability analysis \cite{shufrin2009elastic}, investigations on the sun core temperatures \cite{chandrasekhar1957introduction}, and others, the Bratu equation stands as a significant model in various scientific domains. However, the nonlinear nature of the Bratu equation presents significant challenges for accurate and efficient numerical methods, especially for obtaining sharp solutions. Furthermore, going to higher dimensions increases the difficulty of accurately identifying such sharp solutions.
	
	%\subsection{Numerical Methods 3D Bratu}
	Most research has focused on the 1D or 2D Bratu equation, while studies on the 3D or higher-dimensional cases remain very limited. To our knowledge, only six studies have explored the 3D case, with none extending beyond this dimension.
	% McGough
	The first study conducted by McGough in 1998 \cite{mcgough1998numerical} employed a finite difference method and discretized the domain into a grid of $32^3$ points (including boundary points). McGough then solved the resulting nonlinear system using Newton's method and applied the generalized minimal residual method (GMRES) \cite{saad1986gmres} to address the linear system arising from the inversion of the Jacobian matrix.
	
	% Liao
	In 2012, Liao \cite{liao2012homotopy} demonstrated that the homotopy analysis method (HAM) could be used to transform the 3D Bratu equation into an infinite number of 6th-order linear PDEs which are easily solvable through algebraic calculations. This approach enabled an approximation of the parameter as a function of the infinity norm of the solution.
	
	% Karkowski
	In 2013, Karkowski \cite{karkowski2013numerical} explored three different numerical methods: the pseudospectral method, the finite difference method, and the radial basis functions method. Based on his experiments, Karkowski found the pseudospectral method to be more effective in three dimensions. He implemented this method using $25^3$ grid points chosen from the Chebyshev points of the second kind. 
	The second derivatives at these points were approximated using the second derivative of the Lagrange polynomial interpolant. This resulted in a sparse asymmetric system of linear equations, which can be solved using either the WSMP \cite{gupta2000wsmp} or PARDISO \cite{schenk2004solving}.
	
	% Hajipour
	In 2018, Hajipour et al. \cite{hajipour2018accurate} proposed a fourth-order nonstandard
	compact finite difference formula to discretize the second-order derivative,
	thereby transforming the Bratu equation into a nonlinear algebraic system.
	They developed an iterative method based on Newton's method to solve this
	system. They divided the domain into $15^{3}$ cells and considered
	two parameter values: 0.5 and 1. The infinity norms for the lower solutions
	were found to be 0.02867625 and 0.05855645, while the upper solutions yielded
	norms of 10.58986233 and 8.67456449.
	
	% Iqbal
	In 2020, Iqbal et al. \cite{iqbal2020numerical} utilized the full approximation storage (FAS), a hybrid method combining finite difference and multigrid approaches. The system of equations was solved at each level of the 
	grid throughout the multigrid cycle. To solve the system, they employed the minimal residual method (MINRES) \cite{paige1975solution} as a smoother. Using multigrid, they could utilize up to $161^{3}$ grid points.
	
	% Temimi
	Most recently, in 2024, Temimi et al. \cite{temimi2024efficient} introduced an
	iterative finite difference method. Unlike the standard finite difference
	approach which typically employs Newton's method, their method used
	an iterative scheme to generate a series of solutions. At each step, the
	current solution was fixed for the nonlinear term, and it was then used
	to compute the next solution by solving the resulting linear system. Their
	experiments used a grid of $21^{3}$ points.
	
	%\subsection{Note for all previous work}
	\begin{figure}
		\centering
		\includegraphics[width=0.5\textwidth]{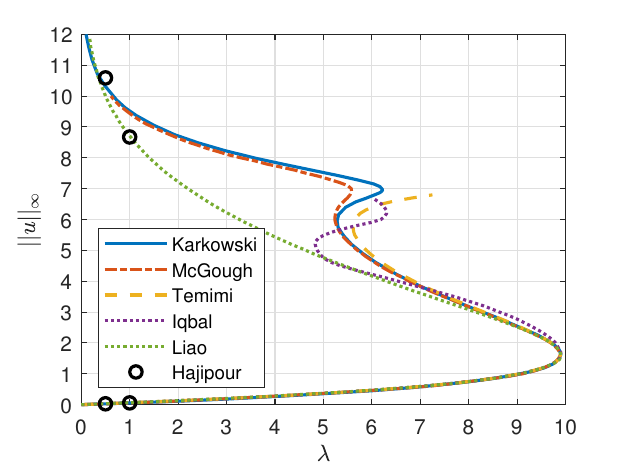}
		\caption{Bifurcation diagrams of the 3D Bratu equation from previous research.}
		\label{fig_3d_bif_prev_research1}
	\end{figure}
	
	We collected the results from these six studies and presented their bifurcation diagrams in Fig.\ \ref{fig_3d_bif_prev_research1}, except Hajipour et al.'s study \cite{hajipour2018accurate}, which provided only four points. Among the methods used to construct the bifurcation diagrams, McGough \cite{mcgough1998numerical} employed arclength continuation, while the others repeatedly varied the parameter value. Notably, the lower solutions, ranging from 0 to the first turning point, appear consistent across all bifurcation diagrams. Note that we need to multiply the parameter of the Bratu equation by 4 for some research that used $[-1,1]^3$ as their domain. Further details can be found in Subsection \ref{connection_domain}.
	
	% Liao and Hajipour
	Unlike the other methods, Liao's approach \cite{liao2012homotopy} did not require discretization or grid points. However, it was only able to identify one turning point. This limitation suggests that the method may lack accuracy, particularly for the upper solutions (i.e., beyond the first turning point), as it fails to capture the second or third turning points. Additionally, Hajipour et al.'s results \cite{hajipour2018accurate}, obtained using a nonstandard compact finite difference formula, show a closer resemblance to Liao's bifurcation.
	
	% 4 lainnya
	Among the remaining studies, Iqbal et al. \cite{iqbal2020numerical} used the highest number of grid points. However, their results did not converge to the correct values as the results displayed different upper solutions or bifurcation diagrams compared to McGough \cite{mcgough1998numerical}, Karkowski \cite{karkowski2013numerical}, and Temimi et al. \cite{temimi2024efficient}. 
	While these three studies showed similar results up to a certain point, they diverged beyond that point due to differences in discretization techniques or methods.
	
	% Temimi
	During our research, we also successfully implemented the iterative scheme proposed by Temimi et al. \cite{temimi2024efficient}, which produced identical results to those obtained by the standard finite difference method solved using Newton's method. However, we observed that Temimi et al.'s iterative method 
	took more iterations to converge. From our perspective, their numerical
	simulation was the simplest compared to McGough and Karkowski's, as they
	used only $21^{3}$ grid points. Moreover, they employed fewer grid points
	than McGough and did not continue to compute the third turning point as
	McGough did.
	
	%\subsection{Intro berikutnya}
	Overall, McGough \cite{mcgough1998numerical} and Karkowski \cite{karkowski2013numerical} appear to have produced better bifurcation diagrams compared to the others. However, the accuracy of their results is limited due to the incorporation of a small number of grid points. Upper solutions with very sharp gradients worsen this limitation.
	
	% Tujuan paper ini
	This paper aims to provide improved solutions and bifurcation diagrams for the 3D Bratu equation. In addition, we will also provide bifurcation diagrams for the 1D, 2D, 4D, and 5D cases. To achieve this, we propose using a finite difference method that incorporates the symmetry properties to reduce its complexity and computational time required for the search process. Additionally, we will numerically investigate the relationship between the Bratu equation on the cube and ball domains.
	
	% Manfaat symmetry dan PD yang memiliki symmetry
	We focus on the symmetry properties in the solutions to the Bratu equation \cite{boyd1986analytical}, the central equation discussed in this paper. Symmetry also appears in many other PDEs. For instance, in the nonlinear Schrödinger equation (NLSE), symmetric solutions such as the Akhmediev breather, the Peregrine breather, and the Kuznetsov-Ma breather are well-known \cite{onorato2013rogue}. Other PDEs, such as the 1D and 2D discrete Allen-Cahn equation with cubic and quintic nonlinearity, also exhibit symmetric steady-state solutions \cite{taylor2010snaking, kusdiantara2019snakes}.
	
	%\subsection{Susunan paper}
	The remainder of this paper is organized as follows:  
	Section \ref{sec2} provides a brief overview of the Bratu equation.  
	Section \ref{sec3} explains the standard finite difference method and its application to different problems.  
	Section \ref{sec4} outlines the proposed method, including how symmetry properties are incorporated.  
	Section \ref{SecResults} presents the numerical results across various dimensions, including solution profiles and bifurcation diagrams.  
	Finally, Section \ref{sec6} summarizes our findings' significance and suggests future research directions.

	\section{Bratu Equation} \label{sec2}
	We consider the Bratu equation \cite{bratu1914equations} in $d$ dimensions, also known as the Gelfand equation, which is a nonlinear elliptic PDE with a positive parameter $\lambda$ and a Dirichlet boundary condition given by
	\begin{equation}
		\begin{split}
			\Delta u(\textbf{x}) + \lambda e^{u(\textbf{x})} & = 0, \qquad \textbf{x}\in\Omega, \\
			u(\textbf{x}) & = 0, \qquad \textbf{x}\in\partial\Omega,
		\end{split}\label{gov}
	\end{equation}
	where $\textbf{x} = (x_1,\dots,x_d)$, $\Omega=\bracketsquare{0,1}^d$, and $\Delta$ is the $d$-dimensional Laplace operator
	\begin{equation}
		\Delta = \sum_{i=1}^{d} \frac{\partial^2}{\partial x_i^2}.
	\end{equation}
	The 1D Bratu equation has an analytical solution given by \cite{shahab2024neural}
	\begin{equation}
		u(x) = 2\ln \bracketround{ \frac{\cosh \theta}{\cosh(\theta(1-2x))} } ,
	\end{equation}
	where $\theta$ is the solution of
	\begin{equation} \label{cosh_theta}
		\cosh \theta = \frac{4\theta}{\sqrt{2\lambda}}.
	\end{equation}
	However, no analytical solutions have been found for higher dimensions.
	
	% Check again dari Temimi, cukup lengkap
	The Bratu equation is solvable for $0<\lambda<\lambda^*$ where $\lambda^*$ represents the first turning point or critical point. The first turning point for the 1D case is 3.513830719 \cite{doedel1997auto97, fedoseyev2000continuation}.
	The first turning points for the 2D case have been computed numerically in several research:
	$1.702 \times 4=6.808$ \cite{boyd1986analytical, liao2012homotopy} and 
	6.808124423 \cite{doedel2000collocation}.
	Similarly, the computed first turning points for the 3D case are
	$2.475 \times 4 = 9.9$ \cite{karkowski2013numerical}, 
	9.90188543 \cite{temimi2024efficient},
	9.90257408 \cite{iqbal2020numerical}, and
	$2.476 \times 4 = 9.904$ \cite{liao2012homotopy}.
	In the next subsection, we will provide insights on obtaining an analytical upper bound for $\lambda$.
	Additionally, we present the approximate first turning points obtained using our proposed method in Section \ref{SecResults} and \ref{appendix_b}.
	
	For the 1D and 2D cases, the Bratu equation has two solutions (lower and upper solutions) for $\lambda\in(0,\lambda^*)$, one solution for $\lambda=\lambda^*$, and no solutions for $\lambda>\lambda^*$. For the 3D case, previous research indicates that it can have more than two solutions for $\lambda\in(0,\lambda^*)$ as shown in Fig.\ \ref{fig_3d_bif_prev_research1}.

	\subsection{Upper Bound for $\lambda$} \label{subsection_upper_bound}
	This section provides an analytical method to obtain an upper bound for $\lambda$. The Bratu equation does not admit a solution if $\lambda$ exceeds this upper bound.
	
	Let $v_1$ be the eigenfunction associated with the first eigenvalue $\lambda_1$ in the spectral problem
	\begin{align} \label{eigenfunctionproblem}
		\begin{split}
			- \Delta v & = \lambda v, \qquad \textbf{x}\in\Omega, \\
			v & = 0, \qquad \;\; \textbf{x}\in\partial\Omega.
		\end{split}
	\end{align}
	Since the first eigenfunction is of one sign, it can be assumed to be positive. Let us normalize it such that
	\begin{equation} \label{normalized_v}
		\int_{\Omega} v_1 \, d\textbf{x} = 1.
	\end{equation}
	Multiplying $- \Delta u = \lambda e^u$ by $v_1$ and integrating over $\Omega$, we obtain
	\begin{equation} \label{eq7}
		\int_{\Omega} - \Delta u \, v_1 \, d\textbf{x} 
		= \lambda \int_{\Omega} e^u v_1 \, d\textbf{x}.
	\end{equation}
	Since both $u$ and $v_1$ are zero on the boundary, we have
	\begin{equation} \label{eq8}
		\int_{\Omega} - \Delta u \, v_1 \, d\textbf{x} 
		= - \int_{\Omega} u \, \Delta v_1 \, d\textbf{x}
		= \lambda_1 \int_{\Omega} u \, v_1 \, d\textbf{x}.
	\end{equation}
	
	Applying Jensen's inequality which states that if $\phi$ is a convex function and $\text{meas}(\Omega) = 1$, then
	\begin{equation}
		\int_{\Omega} v_1 \phi(u) \, d\textbf{x} 
		\ge \phi \bracketround{\int_{\Omega} v_1 u \, d\textbf{x}},
	\end{equation}
	as the norm of $v_1$ is 1. Substituting the exponential function into the inequality, we get
	\begin{equation} \label{eq10}
		\int_{\Omega} e^u \, v_1 \, d\textbf{x}
		\ge e^{ \int_{\Omega} u \, v_1 \, d\textbf{x} }.
	\end{equation}
	From Eqs.\ \eqref{eq7}, \eqref{eq8}, and \eqref{eq10}, we derive the following inequality
	\begin{equation}
		\lambda_1 \int_\Omega u \, v_1 \, d\textbf{x}
		= \lambda \int_{\Omega} e^u v_1 \, d\textbf{x}
		\ge \lambda e^{ \int_\Omega u \, v_1 \, d\textbf{x} }.
	\end{equation}
	Since $e^{f(x)} \ge ef(x)$ \cite{jacobsen2002liouville} for any function $f(x)$, it follows that
	\begin{equation}
		\lambda_1 \int_\Omega u \, v_1 \, d\textbf{x}
		\ge \lambda e \int_\Omega u \, v_1 \, d\textbf{x}.
	\end{equation}
	Thus, we conclude
	\begin{equation}
		(\lambda e - \lambda_1) \int_\Omega u \, v_1 \, d\textbf{x} \le 0.
	\end{equation}
	However, since the solution of the Bratu equation is positive ($u > 0$), a contradiction arises if $(\lambda e - \lambda_1) > 0$. Therefore, we can conclude that the Bratu equation does not admit a solution for $\lambda > \lambda_1/e$.
	
	We can explicitly determine $\lambda_1$ and $v_1$.
	For the eigenvalue problem in $d$ dimension, the eigenfunction is given by
	\begin{equation}
		v = c_1 \sin\bracketround{\sqrt{\frac{\lambda}{d}} \, x_1}
		\cdots \sin\bracketround{\sqrt{\frac{\lambda}{d}} \, x_d}
	\end{equation}
	where $\textbf{x} = (x_1,\dots,x_d)$, and the first eigenfunction occurs when
	\begin{gather}
		\sqrt{\frac{\lambda_1}{d}} = \pi, \qquad \rightarrow \qquad \lambda_1 = d \pi^2.
	\end{gather}
	Hence, the $d$-dimensional Bratu equation does not admit a solution for
	\begin{equation}
		\lambda > \frac{d\pi^2}{e} .
	\end{equation}

	\subsection{A New Constraint $\norm{u}_\infty = A$} \label{sec_new_constraint}
	
	Previous studies show that when the Bratu equation is solved with a numerical method, it usually converges to the lower solution. It is challenging to get the upper solution unless the initial guess is chosen quite close to it.
	
	To address this issue, we found that the original problem \eqref{gov} 
	can be reformulated as follows: given $A$, find $\lambda$ and $u$ that satisfy
	\begin{equation} \label{newBratu}
		\begin{split}
			\Delta u(\textbf{x}) + \lambda e^{u(\textbf{x})} & = 0, \qquad \textbf{x}\in\Omega, \\
			u(\textbf{x}) & = 0, \qquad \textbf{x}\in\partial\Omega, \\
			\norm{u}_{\infty} & = A,
		\end{split}
	\end{equation}
	which will be called by the new Bratu equation.
	The use of $\norm{u}_{\infty} = A$ is based on the fact that multiple solutions exist for each $\lambda \ne \lambda^*$ and each of these solutions can be distinguished by $\norm{u}_{\infty}$. 
	Therefore, fixing $A$ produces a unique solution, making this approach more effective than fixing $\lambda$.
	The advantage of this reformulation is that it allows the upper solution to be obtained directly, starting with initial guesses of $\lambda^{(0)} = 0$ and $u^{(0)} = \textbf{0}$, without requiring guesses close to the upper solution. 
	Additionally, this transformation does not increase the problem's complexity, as the Jacobian matrix's size remains the same. The only difference is that, since $\lambda$ is now a variable in the system, the column representing the derivative of the equations with respect to $\lambda$ in the Jacobian matrix is fully populated.
	Furthermore, for constructing bifurcation diagrams, we will use continuation based on the new Bratu equation by varying the value of $A$ from 0.1 to $A_\text{max}$. This method simplifies the process and reduces the number of iterations required to complete the bifurcation. Although this concept has been used in solving the Bratu Equation on the ball domains \cite{joseph1973quasilinear}, which is an ordinary differential equation explained in the following subsection, to our knowledge, it has not been explored for the higher-dimensional Bratu equation on the cube domains.

	\subsection{Bratu Equation on the Ball Domain}
	The previous discussion focused on the Bratu equation on the cube domain $\Omega=\bracketsquare{0,1}^d$. If we consider the Bratu equation on a $d$-dimensional ball domain  $\Omega = \bracketcurly{ \textbf{x}: \norm{\textbf{x}}_2 \le 1 }$, then all solutions are radially symmetric as stated by Gidas et al. \cite{gidas1979symmetry} and the Bratu equation can be reduced to the following second-order differential equation
	\begin{equation} \label{Bratu_ball_domain}
		u_{\rho\rho} + \frac{d-1}{\rho} u_\rho + \lambda e^u = 0, \qquad 0 < \rho < 1,
	\end{equation}
	with boundary conditions $u_\rho(0)= 0$ and $u(1)=0$, where $\rho$ represents $\norm{\textbf{x}}_2$.
	
	Liouville \cite{liouville1853equation} studied the behavior of the solutions of this equation for $d=1$, while Bratu \cite{bratu1914equations} studied the case for $d=2$. 
	Later, Frank–Kamenetskii \cite{frank1955diffusion}, Chandrasekhar \cite{chandrasekhar1957introduction}, and Gelfand \cite{gelfand1963some} extended the analysis for $d=3$. Finally, Joseph and Lundgren \cite{joseph1973quasilinear} provided the popular complete behavior for all dimensions as follows:
	\begin{enumerate}
		\item For $d=1,2$, there exists $\lambda^* > 0$ such that
		\begin{enumerate}
			\item there are two solutions for $\lambda\in(0,\lambda^*)$,
			\item there is a unique solution for $\lambda=\lambda^*$,
			\item there are no solutions for $\lambda>\lambda^*$.
		\end{enumerate}
		
		\item For $3 \le d \le 9$, let $\tilde{\lambda}=2(d-2)$. Then there exists $\lambda^* > 
		\tilde{\lambda}$ such that
		\begin{enumerate}
			\item there is a finite number of solutions for $\lambda \in (0,\lambda^*) \setminus \{\tilde{\lambda}\}$,
			\item there is a countable infinity of solutions for $\lambda = \tilde{\lambda}$,
			\item there is a unique solution for $\lambda=\lambda^*$,
			\item there are no solutions for $\lambda>\lambda^*$.
		\end{enumerate}
		
		\item For $d \ge 10$, let $\lambda^* = 2(n-2)$. Then
		\begin{enumerate}
			\item there is a unique solution for $\lambda \in (0,\lambda^*)$,
			\item there are no solutions for $\lambda>\lambda^*$.
		\end{enumerate}
	\end{enumerate}
	For $d=2$, it is known that $\lambda^* = 2$.

	\subsection{Connection between Different Domains} \label{connection_domain}
	
	In this subsection, we demonstrate that the value of $\lambda$ from the Bratu equation on $[0,1]^d$ is equivalent to $4\lambda$ for the Bratu equation on $[-1,1]^d$. For simplicity, we will consider the case where $d=1$. Let $u(x)$ be a solution to the equation $u_{xx}+\lambda_1e^u=0$ for $\lambda=\lambda_1$ where $x\in[-1,1]$. Now, define $v(x)=u(2x-1)$ where $x\in[0,1]$. By differentiating, we obtain
	\begin{align}
		%	v_x(x) & = 2u_x(2x-1), \\
		v_{xx}(x) & = 4u_{xx}(2x-1),
	\end{align}
	%or we can write
	%\begin{equation}
	%	u_{xx}(2x-1) = \frac{1}{4}v_{xx}(x).
	%\end{equation}
	and hence
	\begin{align}
		\nonumber
		u_{xx}(2x-1)+\lambda_1e^{u(2x-1)} & =0, \\
		%	\nonumber
		%	\frac{1}{4}v_{xx}(x)+\lambda_1e^{v(x)} & =0, \\
		v_{xx}(x)+4\lambda_1e^{v(x)} & =0,
	\end{align}
	which shows that $v(x)=u(2x-1)$ solves the Bratu equation for $\lambda=4\lambda_1$ where $x\in[0,1]$.

	\section{Finite Difference Method} \label{sec3}
	
	First, it is important to note that the three-point finite difference formulas will be used repeatedly in this research. 
	%Consider $x_i \in [a,b]$, $x_{i+1} = x_i+h$, and $x_{i-1} = x_i-h$ where $h \ne 0$ is sufficiently small to ensure that $x_{i+1}, x_{i-1} \in [a,b]$. Let $u_i$, $u_{i+1}$, and $u_{i-1}$ represent the values of $u$ at $x_i$, $x_{i+1}$, and $x_{i-1}$ respectively. If $u'_i$ and $u''_i$ denote the first and second derivatives of $u$ at $x_i$, then they can be approximated numerically by the three-point finite difference formulas \cite{burden2015numerical}
	%\begin{equation} \label{fd_first}
	%	u'_i \approx \frac{u_{i+1}-u_{i-1}}{2h},
	%\end{equation}
	%and
	%\begin{equation} \label{fd_second}
	%	u''_i \approx 
	%	\frac{u_{i+1} - 2u_i + u_{i-1}}{h^2}.
	%\end{equation}
	Let $u'(x)$ and $u''(x)$ represent the first and second derivatives of $u$ at $x$, respectively. These derivatives can be approximated numerically using the three-point finite difference formulas \cite{burden2015numerical}
	\begin{equation} \label{fd_first}
		u'(x) \approx \frac{u(x+h)-u(x-h)}{2h},
	\end{equation}
	and
	\begin{equation} \label{fd_second}
		u''(x) \approx 
		\frac{u(x+h) - 2u(x) + u(x-h)}{h^2},
	\end{equation}
	for a small value of $h$.
	
	Since this research primarily focuses on higher dimensions, we will use the 3D Bratu equation to illustrate the concept. The explanations provided here can easily be extended to lower or higher dimensions. For simplicity, the term "3D" will be used to refer to both "three-dimensional" and "three dimensions."
	
	The domain for the 3D Bratu equation is $\Omega = [0,1]^3$. Suppose the interval $[0,1]$ is divided into $n$ subintervals of length $h = 1/n$ along the $x$-, $y$-, and $z$-axes. Let $u_{i,j,k}$ and $\Delta u_{i,j,k}$ represent the value and the Laplacian of $u$ at $(x,y,z) = (ih, jh, kh)$ where $i, j, k \in \{0, 1, \dots, n\}$. Extending the finite difference formula for the second derivative in Eq.\ \eqref{fd_second}, the Laplacian in 3D can be approximated by
	\begin{equation}
		\Delta u_{i,j,k} = \frac{u_{i+1,j,k} + u_{i-1,j,k} + u_{i,j+1,k} + u_{i,j-1,k} + u_{i,j,k+1} + u_{i,j,k-1} - 6u_{i,j,k}}{h^2}.
	\end{equation}
	Thus, solving the 3D Bratu equation becomes a problem of finding
	\begin{equation}
		\{u_{i,j,k} \; | \; i,j,k \in \{1,2,\dots,n-1\} \}
	\end{equation}
	that satisfies the system of nonlinear equations
	\begin{equation} \label{nonlinear_system}
		u_{i+1,j,k} + u_{i-1,j,k} + u_{i,j+1,k} + u_{i,j-1,k} + u_{i,j,k+1}
		+ u_{i,j,k-1} - 6u_{i,j,k} + h^2 \lambda e^{u_{i,j,k}} = 0,
	\end{equation}
	for $i,j,k \in \{1,2,\dots,n-1\}$, accompanied with
	\begin{equation} \label{boundary_condition}
		u_{0,j,k}=u_{i,0,k}=u_{i,j,0}=u_{n,j,k}=u_{i,n,k}=u_{i,j,n} = 0, \quad i,j,k \in \{0,1,\dots,n\},
	\end{equation}
	to accommodate the boundary condition. 
	The nonlinear system in Eq.\ \eqref{nonlinear_system} consists of $(n-1)^3$ equations with $(n-1)^3$ variables or unknowns. Then, the problem can be solved with Newton's method or its variants. This approach is generally called a finite difference method (FDM). 
	Note that the Jacobian matrix of this system is sparse, with an average of seven non-zero elements per row or column due to the seven different variables appearing in each equation from Eq.\ \eqref{nonlinear_system}.

	\subsection{Finite Difference Method with $\norm{u}_\infty=A$}
	To solve the new Bratu equation given in Eq.\ \eqref{newBratu} using FDM, we aim to find
	\begin{equation}
		\{\lambda\} \cup \{u_{i,j,k} \; | \; i,j,k \in \{1,2,\dots,n-1\} \} \setminus \{ u_{n/2,n/2,n/2} \}
	\end{equation}
	that satisfies the same nonlinear system as in Eq.\ \eqref{nonlinear_system}, accompanied with 
	\begin{gather}
		u_{0,j,k}=u_{i,0,k}=u_{i,j,0}=u_{n,j,k}=u_{i,n,k}=u_{i,j,n} = 0, \quad i,j,k \in \{0,1,\dots,n\}, \\
		u_{n/2,n/2,n/2} = A.
	\end{gather}
	In this formulation, $\lambda$ is included as one of the unknown variables, while $u_{n/2,n/2,n/2}$ is replaced by the known value from the given condition. This maintains the same number of equations and variables. The condition involving $\norm{u}_\infty$ is replaced by $u_{n/2,n/2,n/2}$ since it is well established that the maximum of $u$ for the Bratu equation always occurs at the center of the domain, which in this case is at $(i,j,k) = (n/2,n/2,n/2)$.
	
	Calculating the Jacobian matrix is the primary difference from the standard Bratu equation. In the new Bratu equation, the column representing the derivatives with respect to $\lambda$ contains full elements. However, this does not alter the fact that the Jacobian matrix remains sparse, with each row having an average of eight non-zero elements.

	\subsection{Finite Difference Method on the Ball Domain}
	
	For the Bratu equation on the ball domain, the interval $\Omega=[0,1]$ is divided into $n$ subintervals of length $h=1/n$. Let $u_i$ represent the value of $u$ at $\rho=ih$, where $i \in {0,1,\dots,n}$. By applying the finite difference formulas for the first derivative in Eq.\ \eqref{fd_first} and the second derivative in Eq.\ \eqref{fd_second}, the problem of solving the Bratu equation on the ball domain is transformed into finding
	\begin{equation}
		\{u_i \; | \; i \in \{1,2,\dots,n-1\} \}
	\end{equation}
	that satisfies the system of nonlinear equations
	\begin{equation} \label{nonlinear_system2}
		\frac{u_{i-1} - 2u_{i} + u_{i+1}}{h^2} + \frac{d-1}{ih} \frac{u_{i+1}-u_{i-1}}{2h} + \lambda e^{u_i} = 0, \quad i \in \{1,\dots,n-1\},
	\end{equation}
	accompanied with
	\begin{equation}
		u_0 = u_1, \quad u_n = 0,
	\end{equation}
	to accommodate the boundary condition. 
	The condition $u_0 = u_1$ ensures that the derivative of $u$ at $\rho=0$ is zero. The nonlinear system in Eq.\ \eqref{nonlinear_system2} consists of $n-1$ equations with $n-1$ variables or unknowns.

	\subsection{Finite Difference Method on the Ball Domain with $\norm{u}_\infty=A$}
	To incorporate the new constraint into the Bratu equation on the ball domain, the problem must be reformulated to find
	\begin{equation} \label{variabel_ball}
		\{\lambda\} \cup \{u_i \; | \; i \in \{2,3,\dots,n-1\} \}
	\end{equation}
	that satisfies the same nonlinear system as in Eq.\ \eqref{nonlinear_system2}, accompanied with 
	\begin{equation} \label{condition_ball}
		u_0 = u_1 = A, \quad u_n = 0.
	\end{equation}
	In this formulation, the term $\norm{u}_\infty$ is replaced by $u_0$ since the maximum of $u$ for the Bratu equation on the ball domain occurs at $i=0$.

	\subsection{Newton's Method} \label{NewtonsMethod}
	In this subsection, we briefly explain how to apply Newton's method in FDM. We illustrate this using the Bratu equation on the ball domain for simplicity.
	
	Let $u = \{u_1, u_2, \dots, u_{n-1} \}^T$ and define the function $F$ as
	\begin{equation}
		F(u) = 
		\begin{pmatrix}
			f_1(u) \\
			f_2(u) \\
			\vdots \\
			f_{n-1}(u)
		\end{pmatrix}
		=
		\begin{pmatrix}
			\frac{u_0 - 2u_1 + u_2}{h^2} + \frac{d-1}{h} \frac{u_2-u_1}{2h} + \lambda e^{u_1} \\
			\frac{u_1 - 2u_2 + u_3}{h^2} + \frac{d-1}{2h} \frac{u_3-u_2}{2h} + \lambda e^{u_2} \\
			\vdots \\
			\frac{u_{n-2} - 2u_{n-1} + u_{n}}{h^2} + \frac{d-1}{(n-1)h} \frac{u_{n}-u_{n-2}}{2h} + \lambda e^{u_{n-1}}
		\end{pmatrix}
	\end{equation}
	which represents the nonlinear system of equations as given in Eq.\ \eqref{nonlinear_system2}.
	The Jacobian matrix of $F$ is defined as
	\begin{equation} \label{jacobian1}
		J(u) =
		\begin{pmatrix}
			\frac{\partial f_1}{\partial u_1} & \frac{\partial f_1}{\partial u_2} & \cdots & \frac{\partial f_1}{\partial u_{n-1}} \\
			\frac{\partial f_2}{\partial u_1} & \frac{\partial f_2}{\partial u_2} & \cdots & \frac{\partial f_2}{\partial u_{n-1}} \\
			\vdots & \vdots & \ddots & \vdots \\
			\frac{\partial f_{n-1}}{\partial u_1} & \frac{\partial f_{n-1}}{\partial u_2} & \cdots & \frac{\partial f_{n-1}}{\partial u_{n-1}} 
		\end{pmatrix} .
	\end{equation}
	Starting from an initial guess $u^{(0)}$, Newton's method proceeds iteratively as follows \cite{chong2023introduction}
	\begin{enumerate}
		\item Solve $J(u^{(k)}) v^{(k)} = F(u^{(k)})$ for $v^{(k)}$,
		\item $u^{(k+1)} = u^{(k)} - v^{(k)}$,
	\end{enumerate}
	for $k=0,1,2,\dots$.
	For this research, the iteration continues until $\norm{u^{(k+1)}-u^{(k)}}_\infty < 10^{-6}$.
	
	For the new Bratu equation on the ball domain where $\norm{u}_\infty = A$, then we set $u = \{\lambda, u_2, u_3, \dots, u_{n-1} \}^T$ and the Jacobian matrix becomes
	\begin{equation}
		J(u) =
		\begin{pmatrix}
			\frac{\partial f_1}{\partial \lambda} & \frac{\partial f_1}{\partial u_2} & \cdots & \frac{\partial f_1}{\partial u_{n-1}} \\
			\frac{\partial f_2}{\partial \lambda} & \frac{\partial f_2}{\partial u_2} & \cdots & \frac{\partial f_2}{\partial u_{n-1}} \\
			\vdots & \vdots & \ddots & \vdots \\
			\frac{\partial f_{n-1}}{\partial \lambda} & \frac{\partial f_{n-1}}{\partial u_2} & \cdots & \frac{\partial f_{n-1}}{\partial u_{n-1}} 
		\end{pmatrix}.
	\end{equation}
	While it is straightforward to construct the sparse tridiagonal Jacobian matrix shown in Eq.\ \eqref{jacobian1} for the one-dimensional case, implementing the sparse Jacobian matrix for higher-dimensional Bratu equations on the cube domains is significantly more complex and challenging.

	\section{Symmetric Finite Difference Method} \label{sec4}
	
	\subsection{Symmetry Property 1}
	
	\begin{figure}
		\centering
		\includegraphics[width=0.5\textwidth]{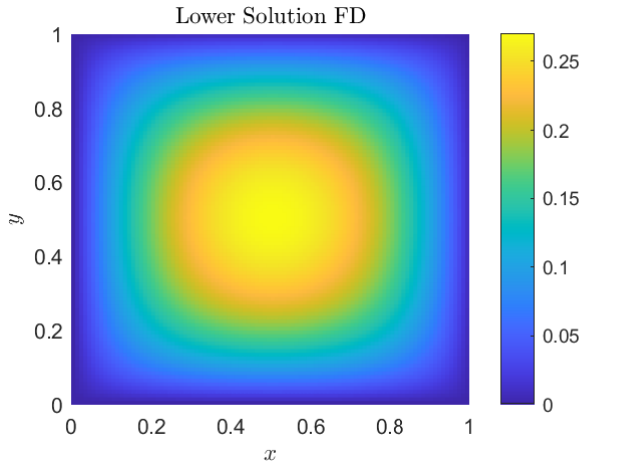}
		\caption{An example of lower solution for the 2D Bratu equation with $\lambda=6$ obtained from FDM with $n=100$.}
		\label{solution_example}
	\end{figure}
	
	Consider the solution of the two-dimensional Bratu equation obtained using FDM as shown in Fig.\ \ref{solution_example}. The solution appears symmetric about both $x=0.5$ and $y=0.5$.
	Therefore, if we know the solution in the region $0 \le x, y \le 0.5$, we can reconstruct the full solution by applying rotations or reflections.
	
	In general, the solution of the $d$-dimensional Bratu equation exhibits symmetry along the planes $x_1=0.5$, $x_2=0.5$, \dots, and $x_d=0.5$. Thus, it suffices to get the solution only for $0 \le x_1,x_2,\dots,x_d \le 0.5$, as the entire solution can be derived from this portion.
	Moreover, dividing the $d$-dimensional cube domain $[0,1]^d$ by the planes $x_1=0.5$, $x_2=0.5$, \dots, and $x_d=0.5$ results in  $2^d$ subdomains. The subdomain corresponding to $0 \le x_1,x_2,\dots,x_d \le 0.5$ is just one of these $2^d$ subdomains.

	\subsection{Symmetry Property 2} \label{symmetric2}
	
	Upon closer inspection of Fig.\ \ref{solution_example}, we observe that the solution is also symmetric about the lines $y=x$ or $y=1-x$. Consequently, the subdomain $0 \le x, y \le 0.5$ can be further divided into two smaller subdomains: $0 \le x \le y \le 0.5$ and $0 \le y \le x \le 0.5$. Let us focus on the subdomain $0 \le x \le y \le 0.5$ for simplicity.
	
	Expanding this concept to the $d$-dimensional cube domain $[0,1]^d$, we uncover a remarkable property: the subdomain $0 \le x_1,x_2,\dots,x_d \le 0.5$ can be represented by an even smaller subdomain $0 \le x_1 \le x_2 \le \cdots \le x_d \le 0.5$.
	
	In the case of the 3D Bratu equation, this property implies that the values of $u$ at a point $(x, y, z)$ and its permutations are equal. Thus, we have the following symmetry
	\begin{equation}
		\begin{split}
			u(x, y, z)
			& = u(x, z, y) \\
			& = u(y, x, z) \\
			& = u(y, z, x) \\
			& = u(z, x, y) \\
			& = u(z, y, x).
		\end{split}
	\end{equation}
	For the $d$-dimensional case, a single value at $(x_1,x_2,\dots,x_d)$ where $0 \le x_1 \le x_2 \le \cdots \le x_d \le 0.5$ can represent $d!$ values in the subdomain $0 \le x_1,x_2,\dots,x_d \le 0.5$. Implementing this concept into FDM with $h=1/n$ for the discretization, which we will refer to as the symmetric finite difference method (SFDM), reduces the number of variables that need to be optimized to approximately
	\begin{equation} \label{number_of_variables0}
		\frac{n^d}{2^d\times d!},
	\end{equation}
	which is significantly smaller than the $(n-1)^d$ variables for the standard FDM without symmetry.
	
	The symmetry properties observed in the numerical solution of the Bratu equation are also inherent in the equation itself. For instance, symmetry is likely present when \cite{bossavit1986symmetry}
	\begin{itemize}
		\item the terms or derivatives with respect to each of the spatial variables are treated equally,
		\item the domain exhibits symmetry, and
		\item the boundary condition also possesses symmetry.
	\end{itemize}
	
	%Last but not least, although some of the presented symmetry properties have been previously studied in \cite{boyd1986analytical, chang2003multigrid, mei1991path}, they only applied the concept to reduce the size of the basis set used to approximate the solutions. Furthermore, none of them tried to incorporate the symmetry properties into FDM like what we will explain in the next parts. In addition, they have not extended into three or higher dimensional problems.
	
	Last but not least, while some of the symmetry properties discussed here have been previously studied in \cite{boyd1986analytical, chang2003multigrid, mei1991path}, those works primarily focused on reducing the size of the basis set used to approximate the solutions. None of them attempted to integrate these symmetry properties into FDM as we will explain in the following parts. Moreover, their studies did not extend to three or higher-dimensional PDEs.
	
	\subsection{Number of Variables} \label{sec_num_of_var}
	In the previous section, we provided an approximate number of variables to be optimized. Here, we present the exact number, particularly for the 3D Bratu equation. The number of variables for other dimensions can be calculated similarly.
	
	Suppose we use $h=1/n$ where $n$ is an even number for discretization. We want to determine the number of indices or points (representing the number of distinct variables) that satisfy $1 \le i \le j \le k \le n/2$. To calculate the number of points, we can divide the problem into four cases.
	The first case is when $(i,j,k)=(a,a,a)$, which gives $n/2$ solutions.  
	The second case is when $(i,j,k)=(a,a,b)$ with $a < b$. There are $n/2$ ways to choose $a$ and $n/2 - 1$ ways to choose $b$. Arranging in increasing order, we get $(n/2)(n/2-1)/2$ solutions in total.  
	The third case is when $(i,j,k)=(a,b,b)$ with $a < b$, yielding $(n/2)(n/2-1)/2$ solutions as well. 
	The fourth case is when $(i,j,k)=(a,b,c)$ with $a < b < c$, which has $(n/2)(n/2-1)(n/2-2)/3!$ solutions.
	Summing these values, we obtain
	\begin{equation}\label{number_of_variables1}
		m = \frac{n^3 + 6n^2 + 8n}{2^3 \cdot 3!} ,
	\end{equation}
	as the number of variables in SFDM for the 3D Bratu equation. This equation is valid for an even number of $n$.
	For an odd number of $n$, we should use $(n-1)/2$, leading to
	\begin{equation} \label{number_of_variables2}
		\frac{(n-1)^3 + 6(n-1)^2 + 8(n-1)}{2^3 \cdot 3!} .
	\end{equation}
	Thus, SFDM with an even number $n=2k$ and an odd number $n=2k+1$ results in the same number of variables to optimize. Furthermore, as $n$ becomes larger, Eqs.\ \eqref{number_of_variables1} and \eqref{number_of_variables2} converge to Eq.\ \eqref{number_of_variables0}.

	\subsection{Symmetric Finite Difference Method}
	
	While explaining the symmetry properties is straightforward, its implementation can be somewhat tricky. The challenge lies in constructing a system of nonlinear equations that incorporate both the symmetry properties and FDM. We will focus on the 3D case to simplify the explanation and provide a step-by-step procedure for constructing and solving this system.

	From the previous explanation, we conclude that by employing the symmetry properties, the 3D Bratu equation can be simplified into a problem of finding
	\begin{equation}
		u = \{u_{i,j,k} \; | \; 1 \le i \le j \le k \le n/2 \}
	\end{equation}
	which satisfies the system of nonlinear equations
	\begin{equation}
		u_{i+1,j,k} + u_{i-1,j,k} + u_{i,j+1,k} + u_{i,j-1,k} + u_{i,j,k+1}
		+ u_{i,j,k-1} - 6u_{i,j,k} + h^2 \lambda e^{u_{i,j,k}} = 0,
	\end{equation}
	for $1 \le i \le j \le k \le n/2$, accompanied with
	\begin{equation}
		u_{0,j,k}=u_{i,0,k}=u_{i,j,0} = 0, \quad i,j,k \in \{0,1,\dots,n/2\}.
	\end{equation}
	However, some cases involve indices outside the set $u$. For example, if $(i,j,k) = (1,1,n/2)$, then $u_{i,j,k+1} = u_{1,1,n/2+1}$ which is undefined and causes the problem to become invalid. This issue must be resolved before proceeding.
	
	Let $m$ denote the number of variables in $u$ that need to be optimized, which can be obtained from Eq.\ \eqref{number_of_variables1} or \eqref{number_of_variables2}. Next, we transform the three-index set $(i,j,k)$, $1 \le i \le j \le k \le n/2$, into a new index using three functions: $f$, $g$, and $M$.
	
	First, the function $f:\mathbb{Z}\rightarrow\mathbb{Z}$ is defined by
	\begin{equation}
		f(i) = \min(i,n-i),
	\end{equation}
	to apply the first symmetry property. Next, the function $g:\mathbb{Z}^3\rightarrow\mathbb{Z}^3$ is defined by
	\begin{equation}
		g(i,j,k) = 
		\begin{cases}
			(i,j,k), & \text{if } i \le j \le k, \\
			(i,k,j), & \text{if } i \le k \le j, \\
			(j,i,k), & \text{if } j \le i \le k, \\
			(j,k,i), & \text{if } j \le k \le i, \\
			(k,i,j), & \text{if } k \le i \le j, \\
			(k,j,i), & \text{if } k \le j \le i, \\
		\end{cases}
	\end{equation}
	to apply the second symmetry property, which sorts the indices in increasing order. Lastly, the function $M:\mathbb{Z}^3\rightarrow\mathbb{Z}$ is defined by
	\begin{equation}
		M(i,j,k) =
		\begin{cases}
			1, & \text{if } (i,j,k) = (1,1,1), \\
			2, & \text{if } (i,j,k) = (1,1,2), \\
			\vdots \\
			\frac{n}{2} + 1, & \text{if } (i,j,k) = (1,2,2), \\
			\frac{n}{2} + 2, & \text{if } (i,j,k) = (1,2,3), \\
			\vdots \\
			m, & \text{if } (i,j,k) = (n/2,n/2,n/2), \\
			0, & \text{if } i=0, \text{or } j = 0, \text{or } k = 0,
		\end{cases}
	\end{equation}
	which maps $(i,j,k)$ to a new index $p$ in increasing order for easier computation and implementation. Additionally, $M$ maps all indices representing boundary points to 0.
	
	Using these functions, we transform the system into
	\begin{multline}
		u_{M(g(f(i+1),j,k))} + u_{M(g(f(i-1),j,k))} + u_{M(g(i,f(j+1),k))} \\ 
		+ u_{M(g(i,f(j-1),k))} + u_{M(g(i,j,f(k+1)))} + u_{M(g(i,j,f(k-1)))} \\
		- 6u_{M(i,j,k)} + h^2 \lambda e^{u_{M(i,j,k)}} = 0,
	\end{multline}
	for $1 \le i \le j \le k \le n/2$, accompanied with
	\begin{equation}
		u_{M(0,j,k)}=u_{M(i,0,k)}=u_{M(i,j,0)} = 0, \quad i,j,k \in \{0,1,\dots,n/2\}.
	\end{equation}
	
	Since there are $m$ equations in total, the system is now reduced to finding
	\begin{equation}
		\{u_1, u_2, \dots, u_m\}
	\end{equation}
	that satisfies the system of nonlinear equations
	\begin{equation} \label{nonlinear_system3}
		-6 u_{p^1_i} + u_{p^2_i} + u_{p^3_i} + u_{p^4_i} + u_{p^5_i} + u_{p^6_i} + u_{p^7_i} + h^2 \lambda e^{u_{p^1_i}} = 0 ,
	\end{equation}
	for $1 \le i \le m$ and $p^1_i, \dots, p^7_i \in \bracketcurly{0, 1, \dots, m}$, accompanied with $u_0=0$ to accommodate the boundary condition.
	
	We can apply Newton's method to solve the system at this point. The contraction of $F$ is straightforward. To speed up the computation of the Jacobian matrix, especially in higher dimensions or with finer grids, we split it into two matrices: one storing the derivatives of the linear terms
	\begin{equation}
		-6 u_{p^1_i} + u_{p^2_i} + u_{p^3_i} + u_{p^4_i} + u_{p^5_i} + u_{p^6_i} + u_{p^7_i} ,
	\end{equation}
	and the other for the nonlinear terms $h^2 \lambda e^{u_{p^1_i}}$. This approach reduces unnecessary computations, as the derivatives of the linear terms remain constant during each iteration of Newton's method (see Appendix A for an example).
	
	Furthermore, we represent the Jacobian matrix as a sparse matrix, which is essential for reducing computational cost when solving $J(u^{(k)}) v^{(k)} = F(u^{(k)})$ for $v^{(k)}$ and for minimizing memory usage. Implementing the sparse matrix representation from the start is crucial, as storing the Jacobian matrix in full before converting it to a sparse format is impractical due to the memory requirements, especially when working with many grid points. Finally, to solve $J(u^{(k)}) v^{(k)} = F(u^{(k)})$ for $v^{(k)}$, we use the "$\backslash$" operator (or the mldivide function) available in Matlab.

	\subsection{Symmetric Finite Difference Method with $\norm{u}_\infty=A$} \label{SFDM_new_Bratu}
	
	To apply SFDM to the new Bratu equation, we need to find
	\begin{equation}
		u = \{\lambda,u_1,u_2,\dots,u_{m-1}\}
	\end{equation}
	that satisfies the same nonlinear system as in Eq.\ \eqref{nonlinear_system3}, accompanied with
	\begin{equation}
		u_0 = 0, \quad u_m = A.
	\end{equation}
	All processes should proceed similarly, except the calculation of the Jacobian matrix, particularly the column representing the derivatives of the equations in the nonlinear system with respect to $\lambda$.

	\section{Experimental Results} \label{SecResults}
	
	Note that all the results presented in the following sections are obtained using Newton's method. The computations were performed on a computer with an Intel Core i5 10400 processor and 64 GB of RAM.
	
	\subsection{Ball Domain}
	
	\begin{figure}
		\centering
		\begin{subfigure}[b]{0.49\textwidth}
			\centering
			\includegraphics[width=\textwidth]{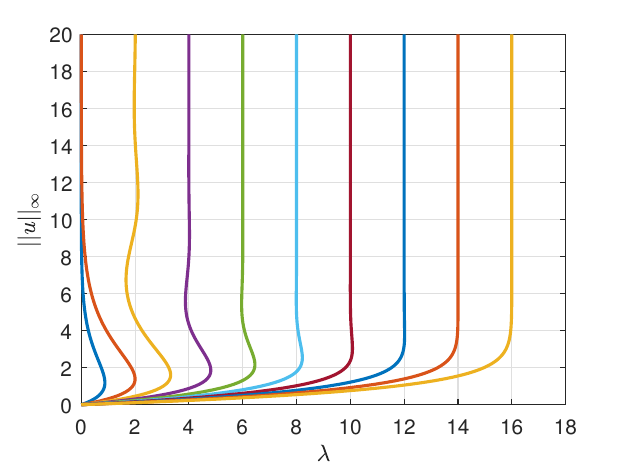}
			\caption{}
			\label{fig_ball_odd_1_sampai_10}
		\end{subfigure}
		\begin{subfigure}[b]{0.49\textwidth}
			\centering
			\includegraphics[width=\textwidth]{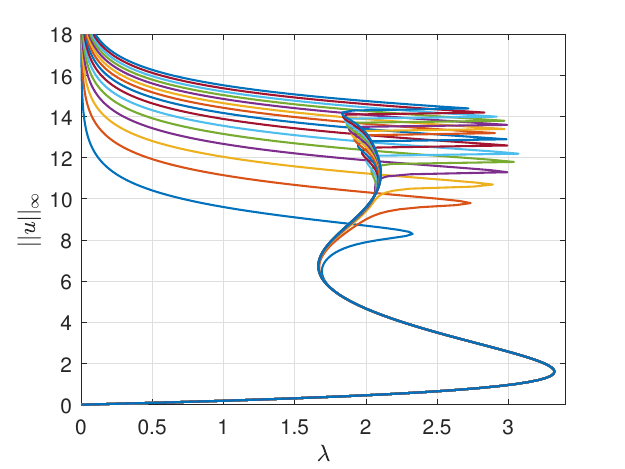}
			\caption{}
			\label{fig_ball_evol}
		\end{subfigure}
		\caption{(a) Bifurcation diagrams of the Bratu equation on the ball domains, from $d=1$ (left) up to $d=10$ (right) with $n=10^6$, (b) bifurcation diagrams for $d=3$ with $n=40,80,\dots,600$.}
		\label{fig_ball_all}
	\end{figure}
	
	We start by presenting the bifurcation diagrams of the Bratu equations on the ball domains, as this provides essential context for the subsequent discussion. In Fig.\ \ref{fig_ball_odd_1_sampai_10}, we display the bifurcation diagrams for dimensions ranging from $d=1$ on the left to $d=10$ on the right. These diagrams were generated by applying continuation on FDM with $n=10^6$ to the new Bratu equation, as described in Eqs.\ \eqref{variabel_ball}, \eqref{nonlinear_system2}, and \eqref{condition_ball}, with $A=0.1, 0.2, \dots, 20$. Notably, all the results are consistent with the behavior summarized by Joseph and Lundgren \cite{joseph1973quasilinear} in the previous section.
	
	We will use these bifurcation diagrams to compare them with those of the Bratu equation on the cube domains. Furthermore, since the first turning points differ between the Bratu equation on the cube and ball domains, we adjust the value of $\lambda$ obtained from the ball domain to align with those of the cube domain. For instance, if the first turning points for the Bratu equation on the ball and cube domains are 2 and 6, respectively, we scale the values of $\lambda$ from the ball domain by a factor of 3.
	
	Although the bifurcation diagrams appear to agree with the findings of Joseph and Lundgren \cite{joseph1973quasilinear}, this alignment is primarily due to using a high number of grid points, $n=10^6$. However, using a smaller value of $n$ produces different results. For example, we generate bifurcation diagrams for $d=3$ with $n=40, 80, \dots, 600$ as shown in Fig.\ \ref{fig_ball_evol}. These bifurcation diagrams do not converge to the expected bifurcation, unlike in Fig.\ \ref{fig_ball_odd_1_sampai_10}. As discussed later, this phenomenon is also observed in the Bratu equation on the cube domains.

	\subsection{3D Bratu Equation}
	
	%\subsubsection{Comparison with Previous Research}
	
	\begin{figure}
		\centering
		\includegraphics[width=0.5\textwidth]{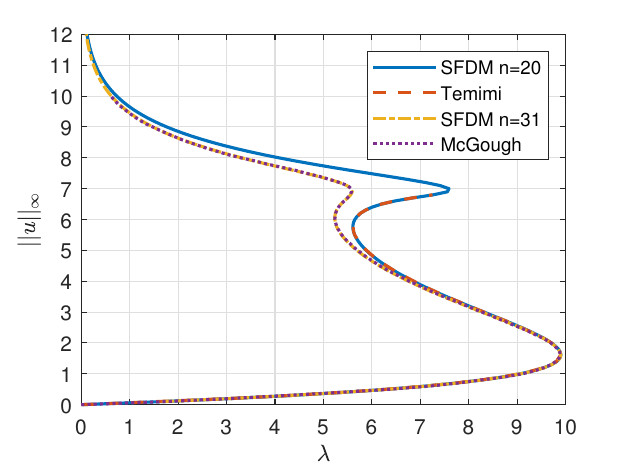}
		\caption{Comparing previous bifurcation diagrams of the 3D Bratu equation with the proposed SFDM with $n=20$ and $n=31$.}
		\label{fig_3d_bif_prev_research_compare}
	\end{figure}
	
	To demonstrate the consistency of our proposed method with the standard FDM, we present the bifurcation diagrams obtained using SFDM with $n=20$ and $n=31$ in Fig.\ \ref{fig_3d_bif_prev_research_compare}. These values of $n$ were selected based on choices made by Temimi et al. \cite{temimi2024efficient} and McGough \cite{mcgough1998numerical} in their research. The comparison reveals that SFDM produces identical bifurcation diagrams, indicating that it yields the same solutions as the standard FDM without employing symmetry. Furthermore, this outcome demonstrates that the iterative approach proposed by Temimi et al. \cite{temimi2024efficient} produces results identical to those obtained by Newton's method.
	
	%\subsubsection{Bifurcation 300}
	
	\begin{figure}
		\centering
		\begin{subfigure}[b]{0.49\textwidth}
			\centering
			\includegraphics[width=\textwidth]{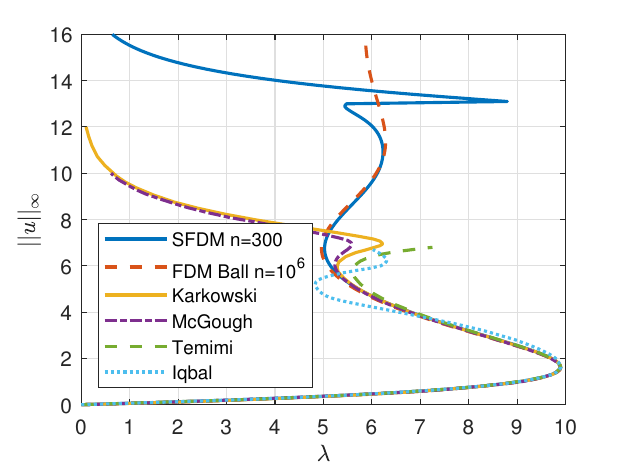}
			\caption{}
			\label{fig_3d_bifurcation_compare_pastresearch}
		\end{subfigure}
		\begin{subfigure}[b]{0.49\textwidth}
			\centering
			\includegraphics[width=\textwidth]{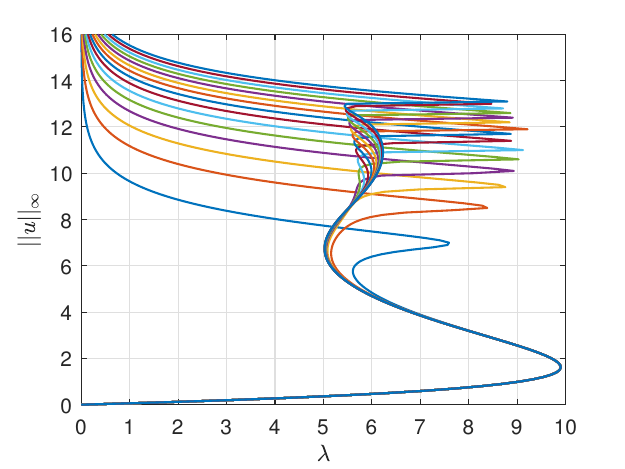}
			\caption{}
			\label{fig_3d_bifurcation_all}
		\end{subfigure}
		
		\caption{(a) Bifurcation diagrams of the 3D Bratu equation from SFDM with $n=300$ compared to the previous research and FDM ball, (b) bifurcation diagrams from SFDM with $n=20, 40, \dots, 300$.}
		\label{fig_3d_bifurcation}
	\end{figure}
	
	The simulation in 3D is limited to a maximum of $n=300$ due to the extensive computational time required. Convergence of SFDM for a single parameter $A$ takes approximately 410 seconds. To construct the bifurcation diagram, we iterate over $A=0.1,0.2,\dots,16$, resulting in a total runtime of around 18 hours. However, SFDM performs significantly faster for smaller values of $n$. For instance, solving for $n=100$ takes less than one second, whereas other methods may struggle with the computational demands of the Bratu equation at these grid points.
	
	The computational time can increase due to a sharp turning point, particularly the last turning point, which poses considerable challenges. Throughout the simulation, we encountered instances where Newton's method failed to converge at this point. In such cases, we replace the previous solution with zeros, $\lambda^{(0)}=0$ and $u^{(0)}=0$. Fortunately, Newton's method can achieve convergence by employing the new Bratu equation. This underscores the effectiveness of the new Bratu equation in constructing bifurcation diagrams.
	
	It is important to note that when the 3D Bratu equation has a similar number of variables as the 1D or 2D cases, the computational time for the 3D case is longer. This is primarily due to the more significant number of elements in the Jacobian matrix, which increases the time required to solve $J(u^{(k)})v^{(k)} = F(U^{(k)})$ for $v^{(k)}$.
	
	In terms of variables, SFDM with $n=300$ has $m=573800$ variables. Compared to the standard FDM with $(n-1)^3=26730899$ variables, SFDM has only around 2.1\% of the variables in the standard FDM.
	
	The complete bifurcation diagram from SFDM with $n=300$ is presented in Fig.\ \ref{fig_3d_bifurcation_compare_pastresearch}. All previously obtained bifurcation diagrams are also included for comparison. A clear difference can be observed between SFDM and the others. 
	The results confirm that the previously obtained solutions were inaccurate, particularly for higher values of $A$. Our findings also reveal new insights into the behavior of the 3D Bratu equation. For example, while previous studies reported only three turning points, we identified five in our experiments. The first turning point occurs at $\lambda^*=9.900146746$, followed by turning points at 5.023545606, 6.229843686, 5.434945323, and 8.797604261.
	Moreover, as shown in Tab.\ \ref{tab:3dtimeratio} in \ref{appendix_b}, the first turning point decreases as $n$ increases. This indicates that the previously reported first turning points of 9.90188543 \cite{temimi2024efficient}, 9.90257408 \cite{iqbal2020numerical}, and 9.904 \cite{liao2012homotopy} are inaccurate, as they are greater than the values we obtained.
	
	Next, we plot the bifurcation diagrams for different grid points, $n=20, 40, \dots, 300$. As $n$ increases, the bifurcation diagrams rise. These bifurcation diagrams show qualitative similarities to those obtained from the ball domain with $d=3$, as illustrated in Fig.\ \ref{fig_ball_evol}. Another critical observation is that as $n$ increases from 20 to 300, the bifurcation diagrams approach those of the Bratu equation on the ball domain with $n=10^6$. It is important to scale the parameter $\lambda$ from the ball domain so that the first turning point of $\lambda$ from the ball domain aligns with that of the cube domain. This demonstrates that as $n$ increases, the qualitative behavior of the Bratu equation on the cube domain becomes similar to that on the ball domain.
	
	%\subsubsection{Accuracy of Bifurcation Diagrams}
	
	\begin{figure}
		\centering
		\begin{subfigure}[b]{0.49\textwidth}
			\centering
			\includegraphics[width=\textwidth]{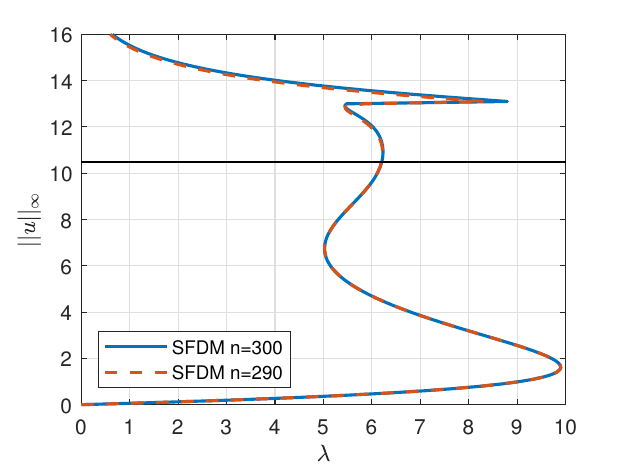}
			\caption{}
			\label{fig_3d_bifurcation_compare290}
		\end{subfigure}
		\begin{subfigure}[b]{0.49\textwidth}
			\centering
			\includegraphics[width=\textwidth]{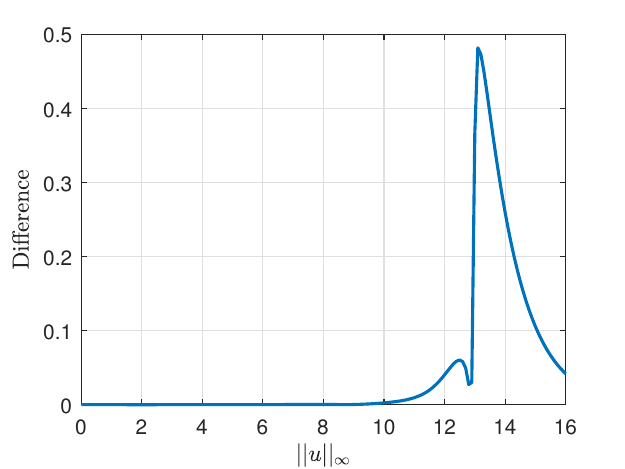}
			\caption{}
			\label{fig_3d_bifurcation_compare290_error}
		\end{subfigure}
		\caption{(a) Bifurcation diagrams of the 3D Bratu equation from SFDM with $n=300$ and $n=290$, the black line is at a height of 10.5, (b) the difference between both bifurcation diagrams.}
		\label{fig_3d_bifurcation_compare290_all}
	\end{figure}
	
	Although we can obtain the bifurcation diagram up to $\norm{u}_\infty = 16$, some upper solutions do not converge to the continuous problem. We compare two consecutive bifurcations to assess the accuracy of the bifurcation diagrams. For example, in the 3D case, we compare results for $n=300$ and $n=290$. We consider the results accurate if the difference between these bifurcation diagrams is less than $0.005$. In other words, for a given value of $A$, the difference in $\lambda$ between $n=300$ and $n=290$ should be less than $0.005$. Using this criterion, we conclude that the bifurcation diagram obtained by SFDM with $n=300$ is accurate up to a height of $10.5$, as shown in Fig.\ \ref{fig_3d_bifurcation_compare290_all}.
	
	%\subsubsection{Benchmark at $\lambda=6$}
	
	\begin{figure}
		\centering
		\begin{subfigure}[b]{0.49\textwidth}
			\centering
			\includegraphics[width=\textwidth]{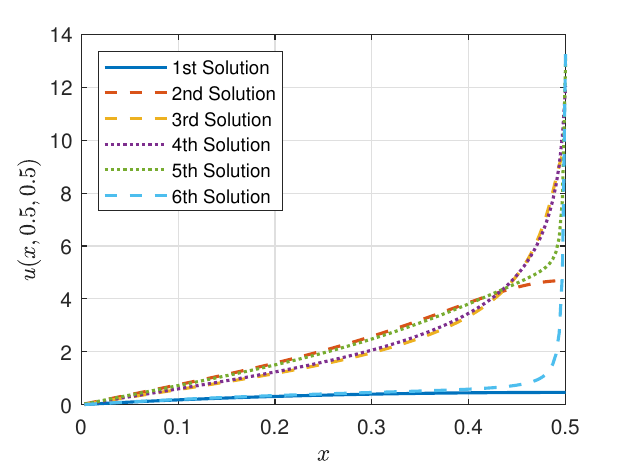}
			\caption{}
			\label{fig_3d_6solutions}
		\end{subfigure}
		\begin{subfigure}[b]{0.49\textwidth}
			\centering
			\includegraphics[width=\textwidth]{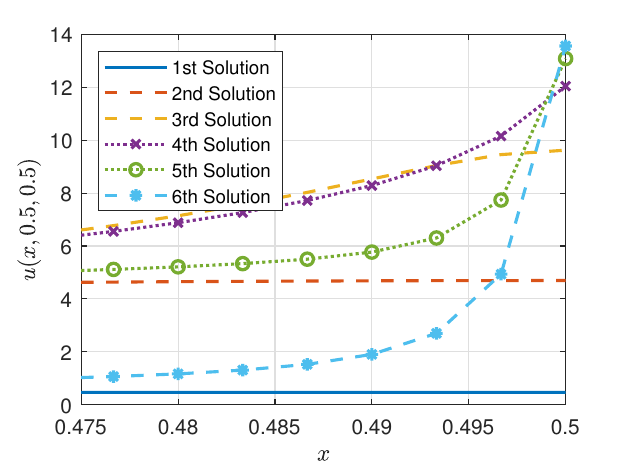}
			\caption{}
			\label{fig_3d_6solutions_zoom}
		\end{subfigure}
		\caption{Six different solutions at $\lambda=6$ for $y=z=0.5$ and (a) $x\in[0,0.5]$, (b) $x\in[0.475,0.5]$.}
		\label{fig_3d_6solutions_all}
	\end{figure}
	
	Next, we present six different solutions obtained from SFDM with $n=300$ for $\lambda=6$, as shown in Fig.\ \ref{fig_3d_6solutions_all}. 
	%While constructing the bifurcation diagram, Newton's method stops when $\norm{u^{(k+1)}-u^{(k)}}_\infty$ is less than $10^{-6}$. For this step, we let it run until the norm is less than $10^{-12}$.
	We plot the solutions for $x\in[0, 0.5]$ and $y=z=0.5$ for easier visualization. The solutions can be reflected across the line $x=0.5$ to obtain the solutions for $x\in[0.5, 1]$. This example shows that the fifth and sixth solutions exhibit very sharp gradients. It is worth noting that sharp gradients occur at all points close to $(x,y,z)=(0.5,0.5,0.5)$ in 3D. This presents a challenge when solving the 3D Bratu equation using methods based on function approximation, such as neural networks \cite{shahab2024neural, raissi2019physics, han2018solving, putri2024deep} and extreme learning machines \cite{fabiani2021numerical}.
	
	%\subsubsection{Even vs Odd}
	
	\begin{figure}
		\centering
		\begin{subfigure}[b]{0.49\textwidth}
			\centering
			\includegraphics[width=\textwidth]{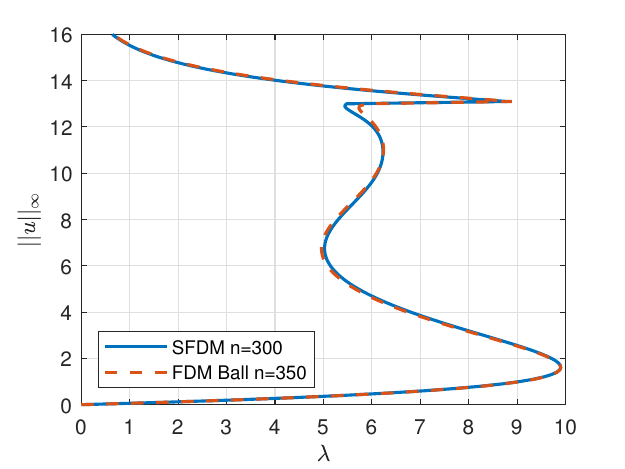}
			\caption{}
			\label{fig_3d_bifurcation_even}
		\end{subfigure}
		\begin{subfigure}[b]{0.49\textwidth}
			\centering
			\includegraphics[width=\textwidth]{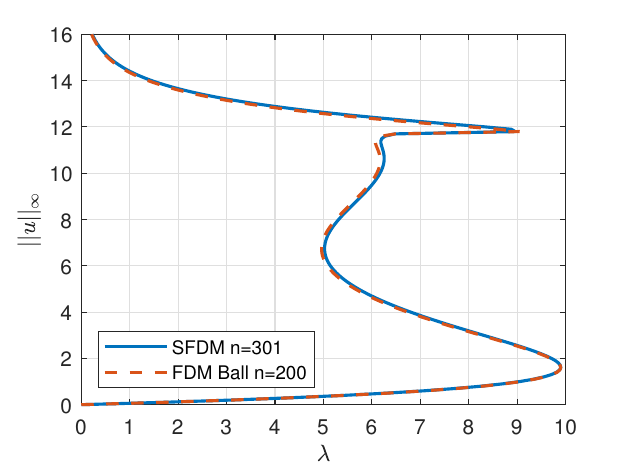}
			\caption{}
			\label{fig_3d_bifurcation_odd}
		\end{subfigure}
		\caption{Bifurcation diagrams of the 3D Bratu equation from SFDM with (a) $n=300$ (even) and (b) $n=301$ (odd). Both are compared to FDM ball with appropriate grids to obtain similar bifurcation diagrams.}
		\label{fig_3d_bifurcation_even_odd}
	\end{figure}
	
	Up to this point, we have consistently used even values for $n$. However, our experiments showed slight differences in the bifurcation diagrams obtained using even and odd values of $n$. To determine the optimal value for $n$, we compare $n=300$ and $n=301$ and plot the resulting bifurcation diagrams. 
	Additionally, we employ FDM for the Bratu equation on the ball domain with $d=3$. 
	Through several trials, we identify values of $n$ on the ball domain that yield similar results to the cube domain: $n=350$ for the even case and $n=200$ for the odd case, as depicted in Fig.\ \ref{fig_3d_bifurcation_even_odd}. This experiment concludes that using an even value for $n$ is more advantageous when solving the Bratu equation, as it provides greater accuracy with the same number of unknowns compared to an odd value of $n$.
	
	%\subsubsection{Stability Analysis}
	
	\begin{figure}
		\centering
		\begin{subfigure}[b]{0.49\textwidth}
			\centering
			\includegraphics[width=\textwidth]{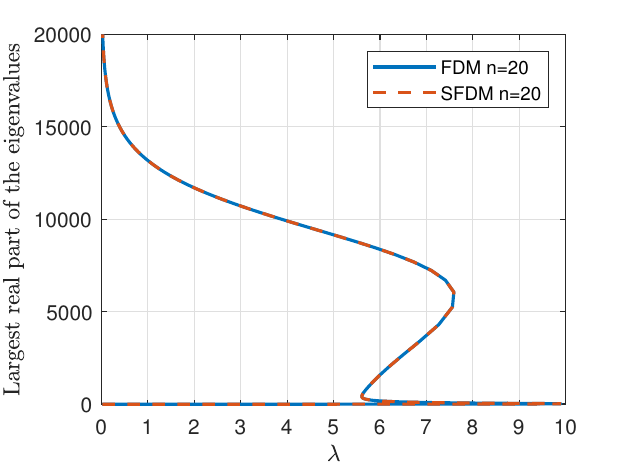}
			\caption{}
			\label{fig_3d_stability_large}
		\end{subfigure}
		\begin{subfigure}[b]{0.49\textwidth}
			\centering
			\includegraphics[width=\textwidth]{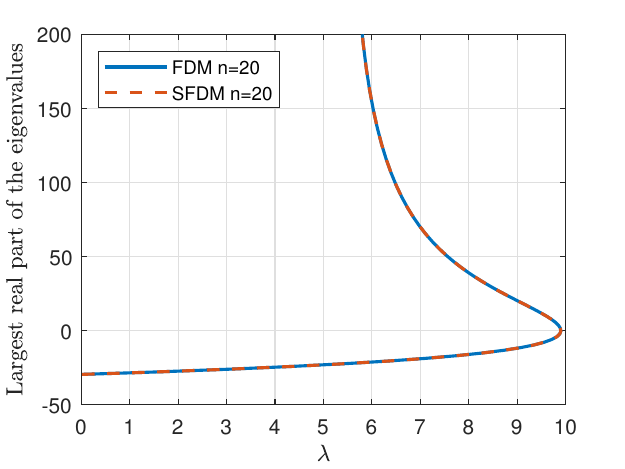}
			\caption{}
			\label{fig_3d_stability_small}
		\end{subfigure}
		\caption{(a) The largest eigenvalue against $\lambda$ from FDM and SFDM with $n=20$, (b) zoom of the graph close to the first turning point.}
		\label{fig_3d_stability}
	\end{figure}
	
	We can delve deeper into the discussion by assuming that the solution of the Bratu equation represents the steady-state solution of $u_t = \Delta u(\textbf{x}) + \lambda e^{u(\textbf{x})}$ and then proceed to determine its linear stability. 
	We choose $n=20$ for this illustration, and the bifurcation diagram is already provided in Fig.\ \ref{fig_3d_bif_prev_research_compare}. To assess the stability of the solutions, we use the Jacobian matrices of the nonlinear systems from both FDM and SFDM. Since these approaches involve different numbers of variables and equations, their Jacobian matrices differ in size. In this case, the Jacobian matrices for FDM and SFDM are $6859\times6859$ and $220\times220$, respectively. We then identify the eigenvalue with the largest real part from these Jacobian matrices. It is important to note that the eigenvalue of the Jacobian matrix is not the same as the parameter $\lambda$ of the Bratu equation.
	
	We illustrate the largest real part of the eigenvalues for both FDM and SFDM in Fig.\ \ref{fig_3d_stability_large}. Remarkably, SFDM yields the same values as those obtained from FDM despite SFDM having a significantly smaller Jacobian matrix. This is intriguing as it shows that even though we reduced the number of variables and equations by modifying FDM into SFDM, the stability properties remain unchanged. Additionally, the computational time required to obtain the SFDM eigenvalues is substantially shorter than for FDM. Moreover, in Fig.\ \ref{fig_3d_stability_small}, the graph intersects the $x$-axis at the first turning point. This suggests that in the 3D case, all the lower solutions are stable, while all the upper solutions are unstable.

	\subsection{4D and 5D Bratu Equation}
	
	\begin{figure}
		\centering
		\begin{subfigure}[b]{0.49\textwidth}
			\centering
			\includegraphics[width=\textwidth]{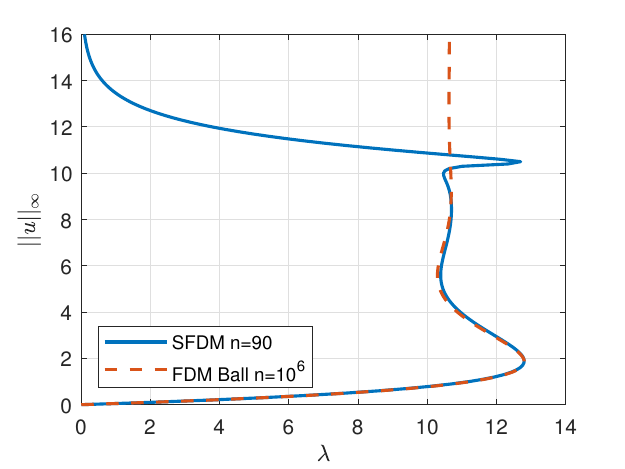}
			\caption{}
			\label{fig_4d_bifurcation}
		\end{subfigure}
		\begin{subfigure}[b]{0.49\textwidth}
			\centering
			\includegraphics[width=\textwidth]{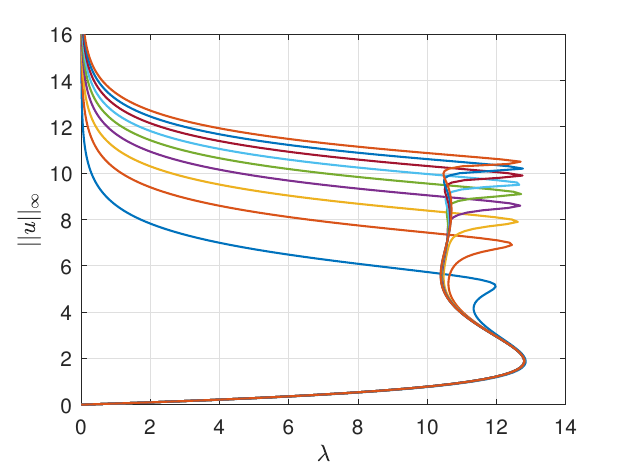}
			\caption{}
			\label{fig_4d_bifurcation_all}
		\end{subfigure}
		\caption{(a) Bifurcation diagrams of the 4D Bratu equation from SFDM with $n=90$ compared to FDM ball, (b) bifurcation diagrams from SFDM with $n=10, 20, \dots, 90$.}
		\label{fig_4d_bifurcation_all_all}
	\end{figure}
	
	\begin{figure}
		\centering
		\begin{subfigure}[b]{0.49\textwidth}
			\centering
			\includegraphics[width=\textwidth]{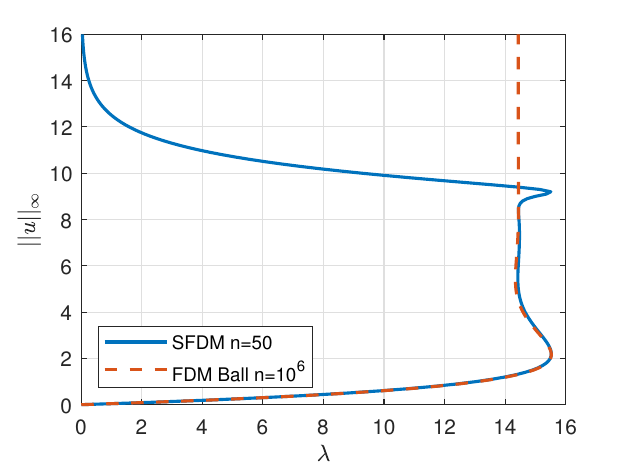}
			\caption{}
			\label{fig_5d_bifurcation}
		\end{subfigure}
		\begin{subfigure}[b]{0.49\textwidth}
			\centering
			\includegraphics[width=\textwidth]{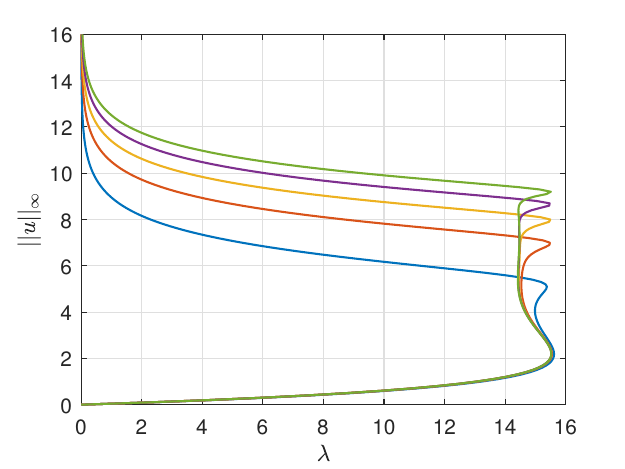}
			\caption{}
			\label{fig_5d_bifurcation_all}
		\end{subfigure}
		\caption{(a) Bifurcation diagrams of the 5D Bratu equation from SFDM with $n=50$ compared to FDM ball, (b) bifurcation diagrams from SFDM with $n=10, 20, \dots, 50$.}
		\label{fig_5d_bifurcation_all_all}
	\end{figure}
	
	Using SFDM, we can further explore the 4D Bratu equation, although the value of $n$ is limited due to the exponential increase in the number of variables. In this study, the maximum tested value for $n$ is 90. Simulating a single value of $A$ takes approximately 380 seconds. To construct the bifurcation diagram, we iterate over $A=0.1, 0.2, \dots, 16$, and generating the complete bifurcation diagram takes about 17 hours. SFDM involves $m=194580$ variables, representing only 0.31\% of the $(n-1)^4=62742241$ variables in the standard FDM.
	
	The bifurcation diagrams of the 4D Bratu equation obtained using SFDM are presented in Fig.\ \ref{fig_4d_bifurcation}, along with a comparison to the bifurcation diagram from the ball domain. Our bifurcation diagram aligns with the ball bifurcation diagram up to a certain height, with the first turning point occurring at $\lambda^*=12.802900147$ (see Tab.\ \ref{tab:4dtimeratio} in \ref{appendix_b}). The evolution of the bifurcation diagrams is illustrated further in Fig.\ \ref{fig_4d_bifurcation_all} for $n=10, 20, \dots, 90$. 
	While in the 3D case, five turning points start to appear when $n=100$, in the 4D Bratu equation, five turning points emerge when $n=50$, but not as distinctly as in the 3D case.
	
	We follow a similar process for the 5D Bratu equation, with $n$ going up to 50. Each simulation for a single value of $A$ takes approximately 350 seconds. We iterate over $A=0.1, 0.2, \dots, 16$ to generate the bifurcation diagram, and completing the entire diagram takes around 16 hours. SFDM involves $m=118755$ variables, which is only 0.042\% of the $(n-1)^5=282475249$ variables from the standard FDM.
	
	The bifurcation diagrams of the 5D Bratu equation obtained by SFDM are shown in Fig.\ \ref{fig_5d_bifurcation}, with a comparison to the bifurcation diagram from the ball domain. Our bifurcation diagram is similar to the ball bifurcation diagram up to a certain height, with the first turning point at $\lambda^*=15.527169368$ (see Tab.\ \ref{tab:5dtimeratio} in \ref{appendix_b}). The evolution of these bifurcation diagrams for $n=10, 20, \dots, 50$ is also depicted in Fig.\ \ref{fig_5d_bifurcation_all}. 
	This is the first research to apply numerical methods to solve the 4D and 5D Bratu equations on the cube domains. From Figs.\ \ref{fig_4d_bifurcation_all_all} and \ref{fig_5d_bifurcation_all_all}, it is evident that as $n$ increases, the bifurcation diagram converges towards the bifurcation diagram of the Bratu equation on the ball domain, similar to the results observed in the 3D case.

	\subsection{1D and 2D Bratu Equation}
	
	\begin{figure}
		\centering
		\begin{subfigure}[b]{0.49\textwidth}
			\centering
			\includegraphics[width=\textwidth]{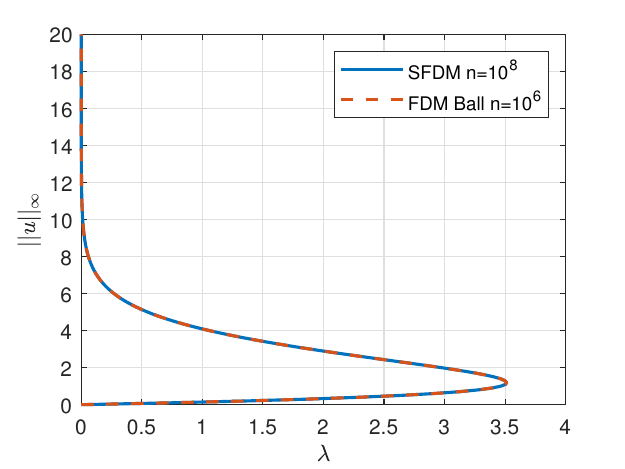}
			\caption{}
			\label{fig_1d_bifurcation}
		\end{subfigure}
		\begin{subfigure}[b]{0.49\textwidth}
			\centering
			\includegraphics[width=\textwidth]{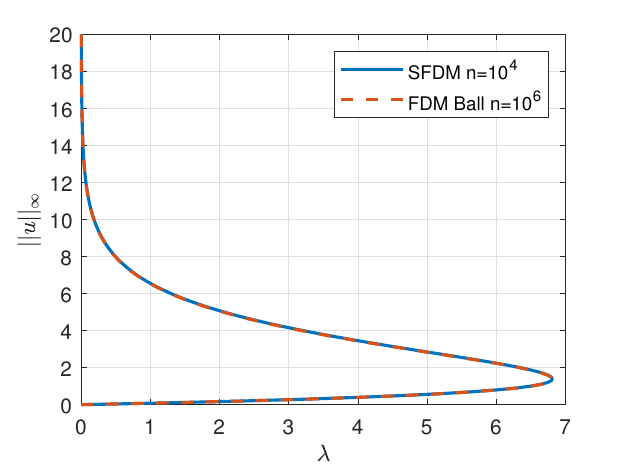}
			\caption{}
			\label{fig_2d_bifurcation}
		\end{subfigure}
		\caption{(a) Bifurcation diagrams of the 1D Bratu equation, (b) bifurcation diagrams of the 2D Bratu equation. Both are compared to FDM ball.}
		\label{fig_1d_2d_bifurcation}
	\end{figure}
	
	We have completed the higher-dimensional part and will shift our focus to the 1D and 2D cases. 
	To solve the 1D Bratu equation, we use SFDM with $n=10^8$. Each computation for a single value of $A$ takes approximately 450 seconds, and constructing the full bifurcation diagram requires about 25 hours. The bifurcation diagram is constructed by iterating over $A=0.1, 0.2, \dots, 20$. In this case, the number of variables is $m=5\times10^7$. The bifurcation diagram is presented in Fig.\ \ref{fig_1d_bifurcation}, along with results from FDM on the ball domain for comparison. 
	The first turning point is obtained as $\lambda^*=3.513830719$ (see Tab.\ \ref{tab14dtimeratio} in \ref{appendix_b}), which matches exactly with the previously reported value in \cite{doedel1997auto97, fedoseyev2000continuation}.
	
	For the 2D case, we employ SFDM with $n=10^4$. Solving for a single parameter value $A$ takes approximately 390 seconds, and completing the bifurcation diagram requires around 22 hours. Like the 1D case, we iterate over $A=0.1, 0.2, \dots, 20$ to build the bifurcation diagram. Here, the number of variables is $m=12502500$, which accounts for about 12.5\% of the $(n-1)^2=99980001$ variables in the standard FDM. The bifurcation diagram is shown in Fig.\ \ref{fig_2d_bifurcation}, alongside results from FDM on the ball domain for comparison. The first turning point is identified as $\lambda^*=6.808124408$. 
	This value matches the previously reported first turning point $6.808124423$ \cite{doedel2000collocation} up to seven decimal places. However, since our first turning point has not yet converged (see Tab.\ \ref{tab:2dtimeratio} in \ref{appendix_b}), and the previous work \cite{doedel2000collocation} used a limited number of mesh and collocation points, it remains unclear which first turning point is more accurate. Nevertheless, we suspect that the exact value of the first turning point lies very close to $6.808124408$ or $6.808124423$.
	
	In the 1D and 2D cases, the bifurcation diagrams on the cube and ball domains are nearly indistinguishable, with only minimal visible differences. We recommend that other researchers use our bifurcation diagrams obtained via SFDM as a reference for future studies on the Bratu equation.

	\subsection{Similar Behavior on Cube and Ball Domains}
	
	Before concluding the discussion, it is important to summarize the general behavior observed in the bifurcation diagrams of the Bratu equation on the cube domains. Based on the numerical results, the following conclusions can be drawn:
	\begin{enumerate}
		\item The bifurcation diagrams for the 1D and 2D Bratu equations on the cube domains are identical to those on the ball domains, as demonstrated in Fig.\ \ref{fig_1d_2d_bifurcation}.
		
		\item In the 3D, 4D, and 5D cases, as $n$ increases, the bifurcation diagrams on the cube domains increasingly resemble those on the ball domains, as illustrated in Figs.\ \ref{fig_3d_bifurcation}, \ref{fig_4d_bifurcation_all_all}, and \ref{fig_5d_bifurcation_all_all}.
		
		\item The evolution of the bifurcation diagrams on the cube domains for small values of $n$ is qualitatively similar to their evolution on the ball domains, as seen in Figs.\ \ref{fig_ball_evol} and \ref{fig_3d_bifurcation_all}.
		
		\item For a given bifurcation diagram on the cube domain with any value of $n=n_1$, it is possible to find a corresponding $n_2$ such that the bifurcation diagram on the ball domain with $n=n_2$ is similar, as shown in Fig.\ \ref{fig_3d_bifurcation_even_odd}.
	\end{enumerate}
	
	While these observations are based solely on the simulations conducted in this research, they provide strong evidence for the following behavior of the Bratu equation on the cube domains:
	\begin{enumerate}
		\item For $d=1,2$, there exists $\lambda^* > 0$ such that
		\begin{enumerate}
			\item there are two solutions for $\lambda\in(0,\lambda^*)$,
			\item there is a unique solution for $\lambda=\lambda^*$,
			\item there are no solutions for $\lambda>\lambda^*$.
		\end{enumerate}
		
		\item For $3 \le d \le 5$, there exists $\tilde{\lambda}$ and $\lambda^* > \tilde{\lambda}$ such that
		\begin{enumerate}
			\item there is a finite number of solutions for $\lambda \in (0,\lambda^*) \setminus \{\tilde{\lambda}\}$,
			\item there is a countable infinity of solutions for $\lambda = \tilde{\lambda}$,
			\item there is a unique solution for $\lambda=\lambda^*$,
			\item there are no solutions for $\lambda>\lambda^*$.
		\end{enumerate}
	\end{enumerate}
	
	The similarity in behavior between the cube and ball domains is particularly intriguing. This parallel, especially in the 3D, 4D, and 5D Bratu equations on the cube domains, has not been previously articulated. However, the exact values of $\tilde{\lambda}$ for each dimension remain unknown. Additionally, we have not explored behavior in higher dimensions due to the lack of simulations. Nevertheless, it is plausible that the behavior in higher dimensions may follow a pattern similar to that observed on the ball domains.

	%\subsection{Accurate First Turning Points}
	%We provide the most accurate first turning points up to this date. We can achieve this because we can use the highest possible grids in all cases of dimensions. To get the first turning points, we do not use interpolations like what many researchers used before. We follow the following procedures to get the first turning points. First, we choose a starting value of $\lambda$, for example we choose 3, 6, and 9 for the starting values in one, two, and three dimensions, respectively. Then we add $0.1$ to the previous value of $\lambda$ to get the new value of $\lambda$, and solve the Bratu equation with Newton's method for the new parameter. If it does not converge, we take the old value of $\lambda$.
	%
	%If it converges, then we add $0.1$ again and solve the Bratu equation with Newton's method again. It it does not converge, 

	\section{Conclusion} \label{sec6}
	
	In this research, we uncover simple yet powerful symmetry properties that can be incorporated into the finite difference method (FDM), which we call the symmetric finite difference method (SFDM). Additionally, we provide a combinatorial approach to calculate the number of variables in SFDM and present a derivation for the new nonlinear system that emerges from it. SFDM effectively reduces the complexity of solving the high-dimensional Bratu equation, enabling a higher number of grid points while decreasing computational time, particularly when used with sparse matrices. Moreover, the introduction of the constraint $\norm{u}_\infty = A$ to the Bratu equation enhances efficiency and accelerates the construction of bifurcation diagrams.
	
	SFDM significantly improves numerical capabilities. In the simulations of the 3D Bratu equation, we utilize up to $301^3$ grid points, surpassing the capabilities of previous research \cite{mcgough1998numerical, liao2012homotopy, karkowski2013numerical, hajipour2018accurate, iqbal2020numerical, temimi2024efficient} and highlighting their inaccuracies, particularly beyond the first turning point. Furthermore, increasing the number of grid points results in bifurcation diagram behavior closely resembling that observed on the ball domains. We also demonstrate that using an even number of $n$ in SFDM yields greater benefits than an odd number. Additionally, simulations show that SFDM achieves similar linear stability analysis despite having fewer variables and equations in the nonlinear system.
	
	Extending SFDM to the 1D, 2D, 4D, and 5D cases further demonstrates its capacity to handle higher grid points and dimensions effectively. The results for these dimensions provide additional evidence that the behavior of bifurcation diagrams for the Bratu equation on the cube domains mirrors the well-known behavior on the ball domains, as documented by Joseph and Lundgren \cite{joseph1973quasilinear}. To our knowledge, this is the first study to emphasize the similarity between the cube and ball domains in the context of the Bratu equation.
	
	The first turning points obtained by SFDM from 1D through 5D are 
	$3.513830719$, $6.808124408$, $9.900146746$, $12.802900147$, and $15.527169368$, respectively. We suggest that future researchers working on the Bratu equation use these turning points as reference values. Lastly, in Appendix B, we provide tables that show the number of variables, computational time, the first turning point, and the ratio of variables in SFDM compared to FDM for various values of $n$ across the 1D to 5D cases.

	%We make the code available online such that it can be used by other research.
	
	In the future, SFDM could be applied to other PDEs that exhibit symmetric solutions similar to those found in the Bratu equation. It would also be interesting to incorporate these symmetry properties into other numerical methods, such as neural networks \cite{shahab2024neural}, spectral methods \cite{boyd1986analytical}, B-spline \cite{caglar2010b}, Green's function \cite{ahmad2023numerical}, or other techniques. Additionally, exploring the properties of other PDEs may reveal further opportunities to enhance numerical methods.

	\section*{CRediT authorship contribution statement}
	\textbf{Muhammad Luthfi Shahab:} Conceptualization, Formal Analysis, Software, Visualization, Writing - original draft, \textbf{Hadi Susanto:} Conceptualization, Supervision, Writing - review \& editing, \textbf{Haralampos Hatzikirou:} Supervision, Writing - review \& editing.

	\section*{Data availability}
	The code and results obtained in this study are available at \\
	\href{https://github.com/luthfishahab/sfdm}{https://github.com/luthfishahab/sfdm}.

	\section*{Declaration of competing interest} 
	The authors declare that they have no known competing financial interests or personal relationships that could have appeared to influence the work reported in this paper.

	\section*{Acknowledgement}
	MLS is supported by a four-year Doctoral Research and Teaching Scholarship (DRTS) from Khalifa University. HS acknowledged support by Khalifa University through a Faculty Start-Up Grant (No.\ 8474000351/FSU-2021-011), a Competitive Internal Research Awards Grant (No.\ 8474000413/CIRA-2021-065), and a Research \& Innovation Grant (No.\ 8474000617/RIG-2023-031). The authors express their gratitude to Prof.\ Mokhtar Kirane for assisting in elucidating the method to obtain the upper bound for $\lambda$, as explained in Subsection \ref{subsection_upper_bound}.

	\bibliographystyle{elsarticle-num} 
	\bibliography{main_references}

	\appendix
	
	\section{Example of SFDM}
	This appendix provides an example of how to derive the nonlinear system of SFDM for the 3D Bratu equation. Let $n=6$, and consequently $m=10$ based on Eq.\ \eqref{number_of_variables1}. For this case, we define $f(i)=\min(i,6-i)$, and the corresponding function $M(i,j,k)$ is given as:
	\begin{equation}
		M(i,j,k) = 
		\begin{cases}
			1, & \text{if } (i,j,k) = (1,1,1), \\
			2, & \text{if } (i,j,k) = (1,1,2), \\
			3, & \text{if } (i,j,k) = (1,1,3), \\
			4, & \text{if } (i,j,k) = (1,2,2), \\
			5, & \text{if } (i,j,k) = (1,2,3), \\
			6, & \text{if } (i,j,k) = (1,3,3), \\
			7, & \text{if } (i,j,k) = (2,2,2), \\
			8, & \text{if } (i,j,k) = (2,2,3), \\
			9, & \text{if } (i,j,k) = (2,3,3), \\
			10, & \text{if } (i,j,k) = (3,3,3), \\
			0, & \text{if } i=0, \text{or } j = 0, \text{or } k = 0.
		\end{cases}
	\end{equation}
	Suppose we take one equation from the system, for example, the equation corresponding to $(i,j,k)=(1,1,3)$ is
	\begin{equation}
		- 6u_{1,1,3} + u_{2,1,3} + u_{0,1,3} + u_{1,2,3} + u_{1,0,3} + u_{1,1,4}
		+ u_{1,1,2} + h^2 \lambda e^{u_{1,1,3}} = 0.
	\end{equation}
	Although in the derivation, we apply $f$, $g$, and $M$ simultaneously, here we demonstrate the step-by-step process for clarity.
	After applying $f$, the equation becomes
	\begin{equation}
		- 6u_{1,1,3} + u_{2,1,3} + u_{0,1,3} + u_{1,2,3} + u_{1,0,3} + u_{1,1,2}
		+ u_{1,1,2} + h^2 \lambda e^{u_{1,1,3}} = 0.
	\end{equation}
	Next, after applying $g$, the equation transforms into
	\begin{equation}
		- 6u_{1,1,3} + u_{1,2,3} + u_{0,1,3} + u_{1,2,3} + u_{0,1,3} + u_{1,1,2}
		+ u_{1,1,2} + h^2 \lambda e^{u_{1,1,3}} = 0.
	\end{equation}
	Finally, after applying $M$, the equation simplifies to
	\begin{equation}
		- 6u_{3} + u_{5} + u_{0} + u_{5} + u_{0} + u_{2}
		+ u_{2} + h^2 \lambda e^{u_{3}} = 0.
	\end{equation}
	Repeating this process for all equations in the system, we obtain
	\begin{equation}
		\begin{split}
			-6 u_{1} + u_{2} + u_{0} + u_{2} + u_{0} + u_{2} + u_{0} + h^2 \lambda e^{u_{1}} & = 0, \\
			-6 u_{2} + u_{4} + u_{0} + u_{4} + u_{0} + u_{3} + u_{1} + h^2 \lambda e^{u_{2}} & = 0, \\
			-6 u_{3} + u_{5} + u_{0} + u_{5} + u_{0} + u_{2} + u_{2} + h^2 \lambda e^{u_{3}} & = 0, \\
			-6 u_{4} + u_{7} + u_{0} + u_{5} + u_{2} + u_{5} + u_{2} + h^2 \lambda e^{u_{4}} & = 0, \\
			-6 u_{5} + u_{8} + u_{0} + u_{6} + u_{3} + u_{4} + u_{4} + h^2 \lambda e^{u_{5}} & = 0, \\
			-6 u_{6} + u_{9} + u_{0} + u_{5} + u_{5} + u_{5} + u_{5} + h^2 \lambda e^{u_{6}} & = 0, \\
			-6 u_{7} + u_{8} + u_{4} + u_{8} + u_{4} + u_{8} + u_{4} + h^2 \lambda e^{u_{7}} & = 0, \\
			-6 u_{8} + u_{9} + u_{5} + u_{9} + u_{5} + u_{7} + u_{7} + h^2 \lambda e^{u_{8}} & = 0, \\
			-6 u_{9} + u_{10} + u_{6} + u_{8} + u_{8} + u_{8} + u_{8} + h^2 \lambda e^{u_{9}} & = 0, \\
			-6 u_{10} + u_{9} + u_{9} + u_{9} + u_{9} + u_{9} + u_{9} + h^2 \lambda e^{u_{10}} & = 0.
		\end{split}
	\end{equation}
	After combining the same terms, substituting $u_0 = 0$, and rearranging, we get
	\begin{equation}
		\begin{split}
			- 6u_{1} + 3u_{2} + h^2 \lambda e^{u_{1}} & = 0, \\
			u_{1} - 6 u_{2} + u_{3} + 2u_{4} + h^2 \lambda e^{u_{2}} & = 0, \\
			2u_{2} - 6u_{3} + 2u_{5} + h^2 \lambda e^{u_{3}} & = 0, \\
			2u_{2} - 6u_{4} + 2u_{5} + u_{7} + h^2 \lambda e^{u_{4}} & = 0, \\
			u_{3} + 2u_{4} - 6u_{5}  + u_{6} + u_{8} + h^2 \lambda e^{u_{5}} & = 0, \\
			4u_{5} - 6u_{6} + u_{9} + h^2 \lambda e^{u_{6}} & = 0, \\
			3u_{4} - 6u_{7} + 3u_{8} + h^2 \lambda e^{u_{7}} & = 0, \\
			2u_{5} + 2u_{7} - 6u_{8} + 2u_{9} + h^2 \lambda e^{u_{8}} & = 0, \\
			u_{6} + 4u_{8} - 6u_{9} + u_{10} + h^2 \lambda e^{u_{9}} & = 0, \\
			6u_{9} - 6u_{10} + h^2 \lambda e^{u_{10}} & = 0.
		\end{split}
	\end{equation}
	
	From this point, we proceed to apply Newton's method. The set of variables to be optimized is $u=\{u_1,u_2,\dots,u_{10}\}$, where
	\begin{equation}
		F(u) = 
		\begin{pmatrix}
			f_1(u) \\
			f_2(u) \\
			\vdots \\
			f_{10}(u)
		\end{pmatrix}
		=
		\begin{pmatrix}
			- 6u_{1} + 3u_{2} + h^2 \lambda e^{u_{1}} \\
			u_{1} - 6 u_{2} + u_{3} + 2u_{4} + h^2 \lambda e^{u_{2}} \\
			\vdots \\
			6u_{9} - 6u_{10} + h^2 \lambda e^{u_{10}}
		\end{pmatrix}.
	\end{equation}
	To expedite the calculation of the Jacobian matrix, particularly for higher grid points, we express the Jacobian as the sum of two matrices: the first matrix stores the derivatives of the linear terms and the second stores those of the nonlinear terms. This results in the following Jacobian:
	\begin{equation}
		J(u) = 
		\begin{pmatrix}
			-6 & 3 & \cdots & 0 \\
			1 & -6 & \cdots & 0 \\
			\vdots & \vdots & \ddots & \vdots \\
			0 & 0 & \cdots & -6 
		\end{pmatrix}
		+
		\begin{pmatrix}
			h^2 \lambda e^{u_{1}} & 0 & \cdots & 0 \\
			0 & h^2 \lambda e^{u_{2}} & \cdots & 0 \\
			\vdots & \vdots & \ddots & \vdots \\
			0 & 0 & \cdots & h^2 \lambda e^{u_{10}} 
		\end{pmatrix}.
	\end{equation}
	In the case of SFDM for the new Bratu equation where $\norm{u}_\infty = A$ as described in Subsection \ref{SFDM_new_Bratu}, the Jacobian matrix takes the following form
	\begin{equation}
		J(u) = 
		\begin{pmatrix}
			0 & -6 & 3 & \cdots & 0 \\
			0 & 1 & -6 & \cdots & 0 \\
			\vdots & \vdots & \vdots & \ddots & \vdots \\
			0 & 0 & 0 & \cdots & 6
		\end{pmatrix}
		+
		\begin{pmatrix}
			h^2 e^{u_{1}} & h^2 \lambda e^{u_{1}} & 0 & \cdots & 0 \\
			h^2 e^{u_{2}} & 0 & h^2 \lambda e^{u_{2}} & \cdots & 0 \\
			\vdots & \vdots & \vdots & \ddots & \vdots \\
			h^2 e^{u_{10}} & 0 & 0 & \cdots & 0
		\end{pmatrix}.
	\end{equation}

	\section{Computational Time, Number of Variables, and First Turning Point} \label{appendix_b}
	
	%This appendix presents tables for the Bratu equation in dimensions ranging from 1D to 5D. These tables illustrate the computational time, the number of variables, and the ratio of variables in SFDM compared to the standard FDM. Specifically, the number of variables in SFDM is denoted as $N_{\text{SFDM}} = m$. In contrast, the number of variables in FDM is given by $N_{\text{FDM}} = (n - 1)^d$, where $n$ represents the discretization parameter. The results are based on solving the new Bratu equation presented in Eq.\ \eqref{newBratu} with $A = 1$. The ratio $r$ is calculated using
	
	%This appendix presents tables for the Bratu equation in dimensions ranging from 1D to 5D. These tables illustrate the number of variables, the computational time, the first turning point, and the ratio of variables in SFDM compared to the standard FDM. The number of variables in SFDM is $N=m$ which can be calculated as in Section \ref{sec_num_of_var}. In contrast, the number of variables in FDM is given by $N = (n - 1)^d$, where $n$ represents the discretization parameter. The computational time is obtained by solving the new Bratu equation presented in Eq.\ \eqref{newBratu} with $A = 1$. The ratio $r$ is calculated using
	
	This appendix presents Tables \ref{tab14dtimeratio}-\ref{tab:5dtimeratio} for the Bratu equation in dimensions ranging from 1D to 5D. The tables show the number of variables, computational time, the first turning point, and the ratio of variables in SFDM compared to the standard FDM. The number of variables in SFDM is $N = m$, which can be calculated as described in Section \ref{sec_num_of_var}. In contrast, the number of variables in FDM is $N = (n - 1)^d$, where $n$ is the discretization parameter. The computational time is determined by solving the new Bratu equation as shown in Eq.\ \eqref{newBratu} with $A = 1$. The ratio $r$ is computed using
	\begin{equation}
		r = \frac{(n-1)^d}{m} .
	\end{equation}
	
	%To get the first turning point we solve the new Bratu equation in 101 different values of $A$ close to the first turning point. For example, in the 1D Bratu equation, the first turning point occurs close to $A=1.2$, so we solve the new Bratu equation for $A=1.1,1.102,\dots,1.3$. After that, we use the Spline interpolation to accurately identifying the first turning point. 
	
	To determine the first turning point, we solve the new Bratu equation for 101 different values of $A$ near the expected first turning point. For instance, in the 1D Bratu equation, the first turning point occurs around $A = 1.2$. Therefore, we solve the equation for values of $A$ ranging from $1.1, 1.102, \dots, 1.3$. After obtaining these solutions, we apply a spline interpolation to accurately identify the exact location of the first turning point.
	
	\newpage
	
	\begin{table}
		\centering
		\caption{The results for the 1D Bratu equation are presented below. The first column lists the values of $n$ used for discretization. The second and third columns show the number of variables in FDM and SFDM, respectively. The fourth and fifth columns provide the computational time (in seconds) and the first turning point obtained using SFDM. The sixth and seventh columns display the ratios $r$ and $1/r$.}
		\label{tab14dtimeratio}
		\begin{tabular}{|c|r|r|r|r|r|r|}
			\hline
			&      FDM &   \multicolumn{3}{c|}{SFDM}    & \multicolumn{2}{c|}{Ratio} \\ \cline{2-7}
			$n$   &      $N$ &      $N$ &  Time & $\lambda^*$ &   $r$ &              $1/r$ \\ \hline
			$10^2$ &       99 &       50 & 0.010 & 3.513647904 & 1.980 &              0.505 \\
			$10^3$ &      999 &      500 & 0.013 & 3.513828891 & 1.998 &              0.501 \\
			$10^4$ &     9999 &     5000 & 0.036 & 3.513830701 & 1.999 &              0.500 \\
			$10^5$ &    99999 &    50000 &  0.33 & 3.513830719 & 1.999 &              0.500 \\
			$10^6$ &   999999 &   500000 &   3.4 & 3.513830719 & 1.999 &              0.500 \\
			$10^7$ &  9999999 &  5000000 &    35 & 3.513830719 & 1.999 &              0.500 \\
			$10^8$ & 99999999 & 50000000 &   450 & 3.513830719 & 1.999 &              0.500 \\ \hline
		\end{tabular}
	\end{table}
	
	\begin{table}
		\centering
		\caption{The same as Table \ref{tab14dtimeratio}, but for the 2D Bratu equation.}
		\label{tab:2dtimeratio}
		\begin{tabular}{|c|r|r|r|r|r|r|}
			\hline
			&      FDM &   \multicolumn{3}{c|}{SFDM}    & \multicolumn{2}{c|}{Ratio} \\ \cline{2-7}
			$n$  &      $N$ &      $N$ &  Time & $\lambda^*$ &   $r$ &              $1/r$ \\ \hline
			100  &     9801 &     1275 & 0.031 & 6.807974209 & 7.687 &              0.130 \\
			200  &    39601 &     5050 & 0.071 & 6.808086880 & 7.842 &              0.128 \\
			400  &   159201 &    20100 &  0.22 & 6.808115038 & 7.920 &              0.126 \\
			500  &   249001 &    31375 &  0.40 & 6.808118416 & 7.936 &              0.126 \\
			1000  &   998001 &   125250 &   1.7 & 6.808122921 & 7.968 &              0.126 \\
			2000  &  3996001 &   500500 &   7.9 & 6.808124047 & 7.984 &              0.125 \\
			4000  & 15992001 &  2001000 &    39 & 6.808124329 & 7.992 &              0.125 \\
			5000  & 24990001 &  3126250 &    67 & 6.808124363 & 7.994 &              0.125 \\
			10000 & 99980001 & 12502500 &   390 & 6.808124408 & 7.997 &              0.125 \\ \hline
		\end{tabular}
	\end{table}
	
	\begin{table}
		\centering
		\caption{The same as Table \ref{tab14dtimeratio}, but for the 3D Bratu equation.}
		\label{tab:3dtimeratio}
		\begin{tabular}{|c|r|r|r|r|r|r|}
			\hline
			&      FDM &   \multicolumn{3}{c|}{SFDM}   & \multicolumn{2}{c|}{Ratio} \\ \cline{2-7}
			$n$ &      $N$ &    $N$ &   Time & $\lambda^*$ &   $r$ &              $1/r$ \\ \hline
			10  &      729 &     35 & 0.0028 & 9.905912320 & 20.83 &              0.048 \\
			20  &     6859 &    220 & 0.0036 & 9.901885432 & 31.18 &              0.032 \\
			30  &    24389 &    680 &  0.011 & 9.900940162 & 35.87 &              0.028 \\
			40  &    59319 &   1540 &  0.027 & 9.900594425 & 38.52 &              0.026 \\
			50  &   117649 &   2925 &  0.067 & 9.900431778 & 40.22 &              0.025 \\
			60  &   205379 &   4960 &   0.11 & 9.900342731 & 41.41 &              0.024 \\
			70  &   328509 &   7770 &   0.20 & 9.900288802 & 42.28 &              0.024 \\
			80  &   493039 &  11480 &   0.34 & 9.900253705 & 42.95 &              0.023 \\
			90  &   704969 &  16215 &   0.58 & 9.900229598 & 43.48 &              0.023 \\
			100 &   970299 &  22100 &   0.87 & 9.900212334 & 43.90 &              0.023 \\
			110 &  1295029 &  29260 &    1.4 & 9.900199548 & 44.26 &              0.023 \\
			120 &  1685159 &  37820 &    2.1 & 9.900189817 & 44.56 &              0.022 \\
			130 &  2146689 &  47905 &    3.4 & 9.900182240 & 44.81 &              0.022 \\
			140 &  2685619 &  59640 &    4.4 & 9.900176226 & 45.03 &              0.022 \\
			150 &  3307949 &  73150 &    7.0 & 9.900171372 & 45.22 &              0.022 \\
			160 &  4019679 &  88560 &    9.7 & 9.900167398 & 45.39 &              0.022 \\
			170 &  4826809 & 105995 &     12 & 9.900164105 & 45.54 &              0.022 \\
			180 &  5735339 & 125580 &     18 & 9.900161344 & 45.67 &              0.022 \\
			190 &  6751269 & 147440 &     26 & 9.900159007 & 45.79 &              0.022 \\
			200 &  7880599 & 171700 &     32 & 9.900157011 & 45.90 &              0.022 \\
			210 &  9129329 & 198485 &     45 & 9.900155294 & 46.00 &              0.022 \\
			220 & 10503459 & 227920 &     55 & 9.900153805 & 46.08 &              0.022 \\
			230 & 12008989 & 260130 &     79 & 9.900152506 & 46.17 &              0.022 \\
			240 & 13651919 & 295240 &    110 & 9.900151366 & 46.24 &              0.022 \\
			250 & 15438249 & 333375 &    120 & 9.900150360 & 46.31 &              0.022 \\
			260 & 17373979 & 374660 &    180 & 9.900149468 & 46.37 &              0.022 \\
			270 & 19465109 & 419220 &    180 & 9.900148673 & 46.43 &              0.022 \\
			280 & 21717639 & 467180 &    230 & 9.900147962 & 46.49 &              0.022 \\
			290 & 24137569 & 518665 &    310 & 9.900147323 & 46.54 &              0.021 \\
			300 & 26730899 & 573800 &    410 & 9.900146746 & 46.59 &              0.021 \\ \hline
		\end{tabular}
	\end{table}
	
	\;
	
	\newpage
	
	\;
	
	\begin{table}
		\centering
		\caption{The same as Table \ref{tab14dtimeratio}, but for the 4D Bratu equation.}
		\label{tab:4dtimeratio}
		\begin{tabular}{|c|r|r|r|r|r|r|}
			\hline
			&      FDM &   \multicolumn{3}{c|}{SFDM}    & \multicolumn{2}{c|}{Ratio} \\ \cline{2-7}
			$n$ &      $N$ &    $N$ &   Time &  $\lambda^*$ &    $r$ &             $1/r$ \\ \hline
			10  &     6561 &     70 & 0.0085 & 12.845620105 &  93.73 &            0.0107 \\
			20  &   130321 &    715 &  0.019 & 12.813772643 & 182.27 &            0.0055 \\
			30  &   707281 &   3060 &   0.10 & 12.807464845 & 231.14 &            0.0043 \\
			40  &  2313441 &   8855 &   0.55 & 12.805226402 & 261.26 &            0.0038 \\
			50  &  5764801 &  20475 &    1.9 & 12.804184914 & 281.55 &            0.0036 \\
			60  & 12117361 &  40920 &     11 & 12.803617732 & 296.12 &            0.0034 \\
			70  & 22667121 &  73815 &     37 & 12.803275250 & 307.08 &            0.0033 \\
			80  & 38950081 & 123410 &    120 & 12.803052770 & 315.62 &            0.0032 \\
			90  & 62742241 & 194580 &    380 & 12.802900147 & 322.45 &            0.0031 \\ \hline
		\end{tabular}
	\end{table}
	
	\begin{table}
		\centering
		\caption{The same as Table \ref{tab14dtimeratio}, but for the 5D Bratu equation.}
		\label{tab:5dtimeratio}
		\begin{tabular}{|c|r|r|r|r|r|r|}
			\hline
			&       FDM &   \multicolumn{3}{c|}{SFDM}    & \multicolumn{2}{c|}{Ratio} \\ \cline{2-7}
			$n$ &       $N$ &    $N$ &   Time &  $\lambda^*$ &     $r$ &            $1/r$ \\ \hline
			10  &     59049 &    126 & 0.0058 & 15.617855802 &  468.64 &          0.00213 \\
			20  &   2476099 &   2002 &  0.068 & 15.547908787 & 1236.81 &          0.00081 \\
			30  &  20511149 &  11628 &    1.2 & 15.534249688 & 1763.94 &          0.00057 \\
			40  &  90224199 &  42504 &     20 & 15.529416008 & 2122.72 &          0.00047 \\
			50  & 282475249 & 118755 &    350 & 15.527169368 & 2378.64 &          0.00042 \\ \hline
		\end{tabular}
	\end{table}

\end{document}